\documentclass[12pt, epsf]{article}
\usepackage{amsfonts, amsmath, latexsym, vmargin, epic, eepic, url}
\usepackage{epsfig, epsf}
\setpapersize{custom}{21cm}{29.7cm}
\setmarginsrb{3cm}{4cm}{3cm}{2cm}{0pt}{0pt}{0pt}{0pt}
\begin{document}

\title{{\bf Data mining for cones of metrics, quasi-metrics, hemi-metrics and super-metrics}}

\author{Michel Deza\footnote{Michel.Deza@ens.fr, CNRS/ENS, Paris and Institute of Statistical Mathematics, Tokyo}~ and Mathieu Dutour \footnote{Mathieu.Dutour@ens.fr, Hebrew University, Jerusalem}\\
}
%,\\
%\'Ecole Normale Sup\'erieure, 45 rue d'Ulm 75230 PARIS Cedex 05, France}

\newtheorem{proposition}{Proposition}
\newtheorem{theor}{Theorem}
\newtheorem{cor}{Corollary}
\newtheorem{lem}{Lemma}
\newtheorem{conj}{Conjecture}
\newtheorem{claim}{Claim}
\newtheorem{remark}{Remark}
\newtheorem{definition}{Definition}
\newcommand{\proof}{\noindent{\bf Proof.}\ \ }

\maketitle

\begin{abstract}
\noindent Using some adaptations of the adjacency decomposition 
method \cite{CR} and the program {\it cdd} (~\cite{Fu}), we compute 
the first computationally difficult cases of convex cones of $m$-ary
and oriented analogs of semi-metrics and cut semi-metrics, which were
introduced in \cite{DR2} and \cite{DP}. We considered also more general 
notion of $(m,s)$-super-metric and corresponding cones. The data on 
related cones - the number of facets, of extreme rays, of their 
orbits and diameters - are collected in Table \ref{tab:MainLovelyTable}. We 
study also criterion of adjacency for skeletons of those cones 
and their duals. Some families of extreme rays and operations on 
them are also given.
\end{abstract}

\section{Introduction}

The notions of $m$-hemi-metric and $m$-partition hemi-metric, as well as of 
quasi-metric and oriented multi-cut quasi-metric, are, respectively, $m$-ary and 
oriented analogs of the (binary and symmetric) notions of metric and cut,
which are important in Graph Theory, Combinatorial Optimization and, more 
generally, Discrete Mathematics.

A finite {\it semi-metric} is a function $d: V_n^2\longrightarrow\mathbb{R}$ with $V_n=\{1, 2, \dots, n\}$ satisfying to $d(x,x)=0$, $d(x,y)=d(y,x)$ and
\begin{equation*}
d(x,y)\leq d(x,z)+d(z,y)\;.
\end{equation*}

Clearly, a semi-metric can take only non-negative values; it is not so,
if we drop the condition $d(x,y)=d(y,x)$.
The notion of a semi-metric was first formalized in the classic paper by 
Fr\'echet \cite{F}. The triangle inequality was first given
as the central property of distances in \cite{F} and later treated in 
Hausdorff \cite{H}. The notion of a metric space appeared also in  
\cite{F}, but the term ``metric'' was first proposed in \cite{H}, page 211.

If the assumption $d(x,y)=d(y,x)$ is replaced by $d(x,y)\geq 0$, then
we obtain 
a {\it quasi-semi-metric}. Those oriented distances appeared already in 
\cite{H}, pages 145--146. Quasi-semi-metrics are used, for example, in the 
Semantics of Computation and in Computational Geometry.

Consider now an extension of the notion of semi-metric in other direction.
A basic example of a metric is $(\mathbb{R}^2, d)$, where $d$ is the Euclidean
distance of $x$ and $y$, i.e. the length of the segment joining $x$
and $y$. An immediate extension is $(\mathbb{R}^3,d)$, where $d(x,y,z)$
is the area of the triangle with vertices $x$, $y$, and $z$.
A $2$-{\it hemi-metric} (see  \cite{M}, \cite{B}, \cite{G1}) is a function $d:V_n^3\longrightarrow\mathbb{R}$ satisfying to
\begin{itemize}
\item[(i)] $d(x,x,y) = 0$, $d(x,y,z)\geq 0$,\\[-8mm]
\item[(ii)] $d(x,y,z)$ is invariant by permutations of $x$, $y$, $z$,\\[-8mm]
\item[(iii)] the $2$-simplex inequality $d(x,y,z)\leq d(x,y,t)+d(x,t,z)+d(t,y,z)$.\\[-8mm]
\end{itemize}
The condition $d(x,y,z)\geq 0$ is convenient for polyhedral methods of the paper; in other contexts some stricter positivity assumptions are necessary.

We recall some terminology. Let $C$ be a polyhedral cone in $\mathbb{R}^{n}$. Given $v \in \mathbb{R}^{n}$, the 
inequality $\sum_{i=1}^nv_ix_i \geq 0$ is said to be {\it valid} for $C$, if it
 holds for all $x \in C$. Then, the set 
$\{ x \in C \mbox{~:~} \sum_{i=1}^nv_ix_i = 0 \}$ 
is called the {\it face of $C$, induced by the valid inequality $\sum_{i=1}^nv_ix_i \geq 
0$}. A face of dimension $\dim(C) - 1$ is called a {\it facet} of $C$;
 a face of dimension $1$ is called an {\it extreme ray} of $C$. An extreme 
ray is called {\em $0,1$ valued} if it contains a vector with only 
$0,1$ values.

Two extreme rays of $C$ are said to be {\it adjacent}, if they
generate a two-dimensional face of
$C$. Two facets of $C$ are said to be {\it adjacent}, 
if their intersection has dimension $\dim(C) - 2$.
The {\it skeleton} graph of $C$ is the graph 
$G(C)$, whose nodes are the extreme rays of $C$ and whose edges
are the pairs of adjacent extreme rays. Denote by $C^*$ the dual cone of $C$. 
The {\it ridge graph} of $C$ is the graph with node set being
the set of facets of $C$ and with an edge between two facets if
they are adjacent on $C$. So, the ridge graph of a cone $C$ is
the skeleton graph $G(C^*)$ of its dual cone. For any cone $C$ we
will call {\em diameter of $C$} the diameter of $G(C)$ and 
{\it diameter of dual $C$} the diameter of its dual, i.e. the
diameter of its ridge graph.

A mapping $f: \mathbb{R}^{n} \longrightarrow \mathbb{R}^{n}$ is called 
a {\it symmetry} of a cone $C$, if it is an isometry, satisfying to
$f(C)=C$ ({\it An isometry} of $\mathbb{R}^{n}$ is a linear mapping 
preserving the Euclidean distance). Given a face $F$, the {\it orbit} 
of $F$ consists of all  faces, that can be obtained from 
$F$ by the group of all symmetries of $C$.
 
For all but one non-oriented cases, the group $Sym(n)$ is {\em the} symmetry
group; but in oriented case appears also a {\it reversal} symmetry 
(see \cite{DDP}); so, all orbits of faces in oriented case are under
action of $Z_2\times Sym(n)$.

The cones $OMCUT_{n}$ and $QMET_{n}$ were introduced and studied, for 
small $n$, in \cite{DP}, \cite{DDP}; their definition and
basic properties are recalled in Section~\ref{TheCaseOfMetricsOrientedOrNot}.

%In cases $SMET^{4, 2}_7$, $HMET_9^6$, 
%$SMET^{4, 3}_7$, $SMET^{5, 5}_8$, $OMCUT_5$, $HCUT^2_6$ and $SCUT^{4, 4}_7$ 

The {\it adjacency decomposition method} of \cite{CR} and the program 
{\it cdd} from \cite{Fu} were used to find the extreme rays and facets
of the hemi-metric and super-metric cones defined in Section 
\ref{BasicDefinitionsOfCones}. When the method succeed completely, we 
indicate this by putting $(a.d.m.)$; note that for the cone 
$HCUT^4_7$, the non-negativity facet was considered as a sub-cone
with its own symmetry group (see \cite{D-PolMet}). If the adjacency 
decomposition method failed to finish (due to the complexity of the
computation), then we get a lower bound on the number of orbits and
we write $\geq$ or $(conj.)$, if, moreover, we expect this lower 
bound to be the right number.
The lower bound for the number of extreme rays of the 
cone $MET_8$, was obtained, using the computation in \cite{DFPS} of the 
vertices of the polytope $MET^{\Box}_8$ (see \cite{DL} page 421). 
Note also that we were able to test adjacency of facets or 
extreme rays of a cone, without knowledge of extreme rays or of 
facets of this cone; namely, we used a linear programming method.

%For the cones $MET_8$, $CUT_8$, and $SMET^{2, 2}_6$ 
%the adjacency decomposition method gave a lower bound, which is 
%expected to be the exact value. 
% For cones 
%$CUT_9$, $MET_9$, $OMCUT_6$, $QMET_6$, $SCUT^{2, 2}_6$, $SMET^{2, 2}_7$, 
%$SMET^{3, 3}_7$, $SMET^{5, 2}_8$, $SMET^{5,3}_8$, $HCUT^5_8$, 
%$SMET^{5,4}_8$ the adjacency decomposition method gives only an 
%estimate of the complexity.

In Table \ref{tab:MainLovelyTable} we summarize the most important numeric 
information on cones under consideration. The column $2$ 
indicates the dimension of the cone, the columns $3$ and $4$ give the number 
of extreme rays and facets, respectively; in parenthesis are given the 
numbers of their orbits.
See \cite{GrMet}, \cite{Gcut} for a description of $MET_7$, $CUT_7$; for all others not bold numbers in the upper half of Table \ref{tab:MainLovelyTable} see \cite{DL}, \cite{DR2}. See \cite{DP} for oriented case with $n\leq 4$.
All data about super-metrics are new; all new data are indicated by bold numbers. The expanded version of those data can be found in \url{http://www.liga.ens.fr/~dutour/}.

\begin{table}
\begin{center}
\scriptsize
\begin{tabular}{|c|c|c|c|c|} 
 \hline \hline
 cone &dim. &Nr. of ext. rays (orbits) &Nr. of facets (orbits)&diameters \\
 \hline
$SMET_{m+2}^{m,s},$&m+2& ${m+2 \choose s+1}$(1)&2m+4(2)&min(s+1,m-s+1); 2\\
 $1 \le s \le m-1$&&&& but 2; 3 if m=2,s=1 \\
 \hline \hline
$SMET_{m+2}^{m,m}$&m+2& m+2(1)&m+2(1)&1; 1 \\
\hline
% $HCUT^m_{m+2}$=$HMET^m_{m+2}$ &m+2& ${m+2} \choose 2$ (1)&$2m+4$ (2)&2; 3,2\\
% for $m \ge 2$&&&&for m=2,>2  \\
% $HCUT^m_{m+3}$&${m+3}\choose 2$& ${m+3}\choose 3$+3${m+3} \choose 4$ (2)&?&?;?\\
% for $m \ge 2$&&&&  \\
% m-$HMET_{m+3}$&${m+3} \choose 2$&?&3${m+3} \choose 2$(2)&?;?\\
% for $m \ge 2$&&&&  \\
 $CUT_4$=$MET_4$ &6&7(2)&12(1)&1; 2 \\
 $CUT_5$ &10&15(2)&40(2)&1; 2 \\
 $MET_5$ &10&25(3)&30(1)&2; 2 \\
 $CUT_6$ &15&31(3)&210(4)&1; 3 \\
 $MET_6$ &15&296(7)&60(1)&2; 2 \\
 $CUT_7$ &21&63(3)&38780(36)&1; {\bf 3} \\
 $MET_7$ &21&55226({\bf 46})&105(1)&3; 2 \\
 $HYP_7$ &21&37170(29)&3773(14)&3; 3 \\
 $CUT_8$ &28&127(4)&$49604520(2169)\;(conj.)$&1; ? \\
 $MET_8$ &28&$119269588({\bf 3918})\;(conj.)$&168(1)&?; 2 \\
 \hline \hline
 $OMCUT_3$=$QMET_3$ &6&12(2)&12(2)&2; 2 \\
 $OMCUT_4$ &12&74(5)&72(4)&{\bf 2; 2} \\
  $QMET_4$ &12&164(10)&36(2)&{\bf 3}; 2 \\
 $OMCUT_5$ &20&540(9)& {\bf 35320(194)} (a.d.m.)&{\bf 2; 3} \\
  $QMET_5$ &20&43590({\bf 229})&80(2)&{\bf 3; 2}\\
 $OMCUT_6$ &30&4682(19)&${\bf \ge 217847040(\ge 163822)}$&{\bf 2}; ?\\
  $QMET_6$ &30&${\bf \ge 492157440(\ge 343577)}$&150(2)&?; {\bf 2}\\
 \hline \hline
 $HCUT_5^2$ &10&25(2)&120(4)&2; {\bf 3} \\
 $HMET_5^2$ &10&37(3)&30(2)&2; 2\\
 $HCUT_6^3$ &15&65(2)&4065({\bf 16})&2; {\bf 3} \\
 $HMET_6^3$ &15&287(5)&45(2)&{\bf 3}; 2  \\
 $HCUT_7^4$ &21&140(2)&{\bf 474390(153)} (a.d.m.)& 2; 3 \\
 $HMET_7^4$ &21&3692(8)&63(2)&{\bf 3}; 2 \\
 $HCUT_8^{5}$ &28&266(2)  &${\bf \ge 409893148(\ge 11274)}$& {\bf 2}; ? \\
 $HMET_8^{5}$ &28&{\bf 55898(13)}  &84(2)& {\bf 3}; {\bf 2} \\
 $HMET_9^{6}$ &36&{\bf 864174(20)} (a.d.m.) &108(2)& ?; 2\\
 $HCUT_6^2$ &20&90(3)&{\bf 2095154(3086)} (a.d.m.)& 2; ? \\
 $HMET_6^2$ &20&12492({\bf 41})&80(2)&{\bf 3}; 2 \\
 $HMET_7^2$ &35&${\bf\ge 454191608(\ge 91836)}$&175(2)&?; 2 \\
 $HMET_7^3$ &35&${\bf\ge 551467967(\ge 110782)}$&140(2)&?; 2 \\
\hline \hline
 $SMET_5^{2,2}$ &10&132(6) &20(1)& 2; 1 \\
 $SCUT_5^{2,2}$ &10&20(2)  &220(6) & 1; 3 \\
 $SMET_6^{3,3/2}$ &15&331989(596)  &45(2)& 6; 2 \\
 $SMET_6^{3,2}$ &15&12670(40)  &45(2)&  4; 2 \\
 $SCUT_6^{3,2}$ &15&247(5)&$866745(1345)\;(conj.)$&  2; ? \\
 $SMET_6^{3,5/2}$ &15&85504(201)  &45(2)& 6; 2 \\
 $SMET_6^{3,3}$ &15&1138(12)  &30(1)& 3; 1 \\
 $SCUT_6^{3,3}$ &15&21(2)  &150(3) & 1; 3 \\
 $SMET_7^{4,2}$ &21&2561166(661) (a.d.m.)&63(2)& ?; 2 \\
 $SMET_7^{4,3}$ &21&838729(274) (a.d.m.)&63(2)& ?; 2 \\
 $SMET_7^{4,4}$ &21&39406(37) &42(1)& 3; 1 \\
 $SCUT_7^{4,4}$ &21&112(2)&148554(114) (a.d.m.)& 1; 4 \\
 $SMET_8^{5,2}$ &28&$\ge 222891598(\ge 6228)$&84(2)& ?; 2 \\
 $SMET_8^{5,3}$ &28&$\ge 881351739(\ge 23722)$&84(2)& ?; 2 \\
 $SMET_8^{5,4}$ &28&$\ge 136793411(\ge 4562)$&84(2)& ?; 2 \\
 $SMET_8^{5,5}$ &28&775807(92) (a.d.m.)&56(1)& ?; 1 \\
 $SMET_9^{6,6}$ &36&$30058078(335)\;(conj.)$&72(1)& ?;1 \\
 $SMET_{10}^{7,7}$ &45&$923072558(1067)\;(conj.)$&90(1)& ?;1\\
 $SMET_6^{2,2}$ &20&$21775425(30827)\;(conj.)$&60(1)& ?; 1 \\
 $SCUT_6^{2,2}$ &20&96(3)&$\ge 243692840(\ge 341551)$& 1; ? \\
 $SMET^{3, 3}_7$ &35&$\ge 594481939(\ge 119732)$&105(1)&?; 1 \\
 $SMET^{2, 2}_7$ &35&$\ge 465468248(\ge 93128)$&140(1)&?; 1 \\
 \hline \hline
\end{tabular}
\caption{Some parameters of cones for small $n$}
\label{tab:MainLovelyTable}
\end{center}
\end{table}

\section{Cuts, semi-metrics and their oriented analogs }\label{TheCaseOfMetricsOrientedOrNot}

We start with the notion of cut semi-metric. Given a subset $S$ of $V_{n}$, let $\delta (S)$ denote the vector in $R^{\frac{n(n-1)}{2}}$, 
defined by $\delta (S)(x,y) = 1$, if $\vert S \cap \{ x, y \} \vert = 1$, 
and $\delta (S)(x,y) = 0$, otherwise (for $1 \leq x < y \leq n$). Obviously, 
$\delta (S)$ defines a semi-metric on $V_{n}$ called a {\it cut semi-metric} (or, simply, a {\it cut}).

Consider now the notion of multi-cut semi-metric. Let $q \geq 2$ be an 
integer and let $S_{1}, \ldots, S_{q}$ be a $q$-partition of $V_{n}$. 
Then the {\it multi-cut semi-metric} $\delta (S_{1}, \ldots, 
S_{q})$ is the vector in  $\mathbb{R}^{\frac{n(n-1)}{2}}$, defined by $\delta (S_{1}, \ldots, 
S_{q})(x,y) = 0$, if $x, y \in S_{h}$ for some $h$, $1 \leq h \leq q$, and 
$\delta (S_{1}, \ldots, S_{q})(x,y) = 1$, otherwise, for $1 \leq x < y \leq 
n$.

It turns out, that every multi-cut is expressed as a sum of cuts with 
non-negative coefficients; so, they are not extreme rays in the 
semi-metric case. The cone defined by all $2^{n-1}-1$ non-zero cuts
is called {\it cut cone} and denoted by $CUT_n$; it 
has full dimension $\frac{n(n-1)}{2}$. The skeleton of $CUT_n$ is 
$K_{2^{n-1}-1}$ (see \cite{DL} pages 539-540 about
$CUT_n$).

The {\it triangle} inequalities $T_{x,y; z}: d(x,y)+d(y,z)-d(x,z)\geq 0$ 
are facets of $CUT_n$ (see \cite{DL}, \cite{GrMet} for details). It is 
conjectured that every facet of $CUT_n$ is adjacent to a triangle 
inequality facet; this conjecture, if true, would imply that the
diameter of dual $CUT_n$ is $3$ or $4$.

The {\em semi-metric cone} $MET_n$ is defined by all triangle inequalities $T_{x,y;z}$. The 
cuts $\delta(S)$ are extreme rays of $MET_n$ and we have the inclusion 
$CUT_n\subset MET_n$, i.e. $MET_n$ is a {\it relaxation} of the cut 
cone $CUT_n$. For $n\leq 4$ both cones coincide, while the number of 
facets of $CUT_n$ and extreme rays of $MET_n$ explodes for $n=7$ 
(see \cite{Gcut} and \cite{GrMet}). 
The diameter of dual $MET_n$ is two (see \cite{DD}), while
the diameter of $MET_n$ (for $n\geq 7$) is three, if the Laurent-Poljak 
conjecture (\cite{LP} and page 543 of \cite{DL}) is true.

%for $HCUT^m_n$ it is $2$ if Conjecture \ref{Conjecture-adjacencies-extrays} holds.

The {\em hypermetric cone} $HYP_n$ is a (smaller than $MET_n$) relaxation of the cut 
cone, which is defined by taking, as the facets, a larger set of faces
of $CUT_n$: the hypermetric inequalities. 
One has $CUT_n\subset HYP_n$ with equality only for $n\leq 6$; $HYP_7$ is the 
first non-trivial instance of this cone (see \cite{DD-hyp} for $HYP_7$ and 
\cite{DL} part II for general $HYP_n$).
\newline

A function $d(x,y)$ is called a {\em quasi-semi-metric} if it satisfies 
the {\it non-negativity} inequalities $NN_{x,y}: d(x,y) \geq 0$ and 
the {\it oriented triangle} inequalities $OT_{x,y; z}: d(x,y) + d(y,z) - d(x,z) \geq 0$.

In the same way, given a subset $S$ of $V_{n}$, let $\delta ^{'} (S)$ denote 
the vector in $R^{n(n-1)}$, defined by $\delta^{'} (S)(x,y) = 1$, if 
$x \in S$, $y \not\in S$, and $\delta^{'} (S)(x,y) = 0$, otherwise, for 
$1 \leq x \not= y \leq n$. Clearly, $\delta^{'} (S)$ defines a 
quasi-semi-metric on $V_{n}$, called an {\it oriented cut}.
Given an ordered $q$-partition $S_{1}, \ldots, S_{q}$ 
of $V_{n}$, let $\delta^{'} (S_{1}, \ldots, S_{q})$ denote the vector in 
$\mathbb{R}^{n(n-1)}$, defined by $\delta^{'} (S_{1}, \ldots, S_{q})_{ij} = 1$, if
 $x \in S_{\alpha}$, $y \in S_{\beta}$, when $\alpha < \beta$, and
 $\delta^{'}(S_{1}, \ldots, S_{q})(x,y) = 0$, otherwise. It may be verified 
that $\delta^{'} (S_{1}, \ldots, S_{q})$ defines a 
quasi-semi-metric on $V_{n}$, which is called an {\it oriented multi-cut}.

The cone generated by all non-zero oriented multi-cuts on $V_n$ is denoted 
by $OMCUT_n$. This cone is of full dimension $n(n-1)$, but, contrary to the 
non-oriented case, oriented multi-cuts are not always adjacent (see \cite{DP} 
and \cite{DDP}), while oriented cuts are still adjacent.
It was conjectured in \cite{DDP} that non-negativity and oriented triangle inequalities are facets of $OMCUT_n$.

Let us denote by $QMET_n$ the cone defined by the $n(n-1)$ non-negativity 
inequalities and the $n(n-1)(n-2)$ oriented triangle inequalities on $V_n$.
It is conjectured in \cite{DDP} that the diameter of $OMCUT_n$ and of dual $QMET_n$ is two.

We have the inclusion $OMCUT_n\subset QMET_n$ with equality if and only 
if $n=3$. The number of facets of $OMCUT_n$ and of extreme rays of $QMET_n$ 
explode for $n=5,6$.

The knowledge of extreme rays of $QMET_n$ 
(the cone of all quasi-semi-metrics on $n$ points) for 
small $n$ will help to build a theory of multi-commodity flows on oriented 
graphs, as well as it was done for non-oriented graphs, using dual $MET_n$.
%(cone of all semi-metrics on $n$ points).

\section{$m$-hemi-metrics, $(m,s)$-super-metrics}\label{BasicDefinitionsOfCones}

The notion of $m$-hemi-metric is an $m$-ary analog of the (binary) 
notion of semi-metric. This notion was introduced in \cite{DR1}
and studied, using the program {\em cdd} \cite{Fu}, for small parameters $m,n$ 
in~\cite{DR2}.

For an arbitrary positive integer $m$, a map $d:
E^{m+1}\longrightarrow\mathbb{R}$ is {\it totally symmetric} if for
all $x_1, \dots, x_{m+1}\in E$ and every permutation $\pi$ of $\{1,\dots, m+1\}$
$$d(x_{\pi(1)}, \dots, x_{\pi(m+1)}) = d(x_1, \dots, x_{m+1}).$$

If it satisfies to $d(x_1, \dots, x_{m+1})=0$, whenever $x_1, \dots, x_{m+1}$ are 
not pair-wisely distinct, then it is called {\it zero-conditioned}.

\begin{definition}
Let $m\geq 1$. An {\it $m$-hemi-metric} on $E$ is a pair 
$(E, d)$, where $d: E^{m+1}\longrightarrow\mathbb{R}$ is totally symmetric,
zero-conditioned and satisfies the {\it $m$-simplex inequality}: for all 
$x_1, \dots, x_{m+2}$$\in E$ 
$$ST_{x_1, \dots,x_{m+2} ; x_{m+2}}: d(x_1, \dots, x_{m+1}) \leq\sum^{m+1}_{i=1}d(x_1, \dots, x_{i-1}, x_{i+1}, \dots, x_{m+2}),$$
and the {\em non-negativity inequality} 
$$NN_{x_1, \dots, x_{m+1}}: d(x_1, \dots, x_{m+1})\geq 0.$$

\end{definition}

\begin{definition}
Let $m$ be a positive integer and let $s$ be any positive number. An {\it $(m,s)$-super-metric} on $E$ is a pair $(E, d)$, 
where $d: E^{m+1}\longrightarrow\mathbb{R}$ is totally symmetric, zero-conditioned and satisfies the {\it $(m,s)$-simplex inequality}:
for all {\em distinct} $x_1, \dots, x_{m+2}\in E$
$$s-ST_{x_1, \dots,x_{m+2} ; x_{m+2}}: sd(x_1, \dots, x_{m+1}) \leq\sum^{m+1}_{i=1}d(x_1, \dots, x_{i-1}, x_{i+1}, \dots x_{m+2}), $$
and the {\em non-negativity inequality} 
$$NN_{x_1, \dots, x_{m+1}}: d(x_1, \dots, x_{m+1})\geq 0.$$
\end{definition}

So, a $m$-hemi-metric is just a $(m,1)$-super-metric and a semi-metric is a $(1,1)$-super-metric.

If $T=\{x_1, \dots, x_{m+2}\}$, then we will set $d_{x_i}=d(x_1,\dots, x_{i-1}, x_{i+1}, \dots, x_{m+2})$ and $\Sigma_{T}=\sum_{i=1}^{m+2} d_{x_i}$. The $(m,s)$-simplex inequality $s-ST_{T; x_{m+2}}$ can be rewritten as $(s+1)d_{x_{m+2}}\leq \Sigma_{T}$.

%One can check also that for $n=m+2$ any partition $(m-s+1)$-hemi-metric,
%such that all, but one, sets $S_{i}$ are singletons, is also
%$(m,s)$-super-metric. 

The following  notation will be used below:
\begin{itemize}
\item the {\it cone $HMET_{n}^m$ of $m$-hemi-metrics}, 
defined by all $(m+2) {n \choose {m+2}}$ \hspace{3mm} $m$-simplex
inequalities and all $n \choose {m+1}$ non-negativity inequalities on $V_n$,
\item the {\it cone $SMET_{n}^{m,s}$ of $(m,s)$-super-metrics},
defined by all $(m+2) {n \choose {m+1}}$ \hspace{3mm}
$(m,s)$-simplex
inequalities and all $n \choose {m+1}$ non-negativity 
inequalities on $V_n$,
\item the {\it cone $SCUT_{n}^{m,s}$}, generated by all $0,1$ valued 
extreme rays of $SMET^{m,s}_n$.
\end{itemize}

All these cones lie in the positive orthant ${\bf R}_+^{n\choose {m+1}}$.

An example of a $m$-hemi-metric is $m$-dimensional volume in $\mathbb{R}^n$ (with $n\geq m$).

The same notion of $m$-volume, restricted to some subset $T$ of 
$\mathbb{R}^n$, gives examples of $(m,s)$-super-metric with 
$s\geq 1$ (see \cite{MaDu}). 
One can check that the $m$-dimensional volume is, moreover, a 
$(m,s)$-super-metric:

(i) if $T$ is $m+2$ vertices of regular $m+1$-simplex and $s=m+1$;

(ii) if $T$ is six vertices of octahedron, $m=2$ and $s=1+\sqrt{3}$.

%(NEED TO CONSIDER ADDING REFERENCE TO MAEHARA)

%(ii) if $T$ is the $2m+2$ vertices of $m+1$-cross-polytope and $s=m-1+\sqrt{m+1}$
%If $d$ is a semi-metric, then $d^2$ is an $(1,1/2)$-super-metric (it follows from $d(x,y)\leq d(x,z)+d(z,y)$ using $2uv\leq u^2+v^2$).
%\begin{equation*}
%d(x,y)\leq d(x,z)+d(z,y)\Rightarrow d^2(x,y)\leq d^2(x,z)+d^2(z,y)+2d(x,z)d(z,y)\leq 2[d^2(x,z)+d^2(z,y)]
%\end{equation*}

Any facet or extreme ray of the cones $SMET^{m,s}_n$ is given by vector, say, $v$, which can be indexed by $(m+1)$-subsets of the set $V_{n}$. So, each such vector can be seen as vertex-labeled subgraph of the Johnson graph $J(n, m+1)$, i.e. we consider the restriction of $J(n,m+1)$ on the {\em support} of $v$ (i.e. the set of all indices of the vector, on which its components are non-zero). Call this restriction, {\em labeled representation graph of $v$} (in Johnson graph $J(n,m+1)$) and denote it by $G_{v}$. In the special case, when $v$ is an $(0,1)$-vector, $G_{v}$ is just usual graph. In Section \ref{case-m-m+3} we will introduce another graph $H_v$ for the special case $n=m+3$.

Below $K_n$ denotes the complete graph on $n$ vertices, $C_n$ denotes the cyclic graph on $n$ vertices, $K_{n_1, \dots, n_t}$ denotes the complete $t$-partite graph on sets of size $n_1,\dots, n_t$. $\nabla G$ is the graph, obtained from a graph $G$ by adding a vertex, which is adjacent to all its vertices, and $K_n-tK_2$ denotes the complete graph $K_n$ with $t$ disjoint edges removed.

\begin{proposition}
For $SMET^{m,s}_n$ holds:

(i) If $s>m+1$, the cone collapses to $0$.

(ii) If $s=m+1$, the cone $SMET^{m,m+1}_n$ collapses to the half-line of all non-negative multiples of the vector of all ones.

(iii) If $m\leq s<m+1$ the non-negativity inequalities are implied by the $(m,s)$-simplex inequalities.
\end{proposition}
\proof Let take the fixed support $T=\{x_1, \dots,x_{m+2}\}$; by summing all $(m,s)$-simplex inequalities $(s+1)d_{x_i}\leq \Sigma_T$, we obtain $(s+1)\Sigma_T\leq (m+2)\Sigma_T$. By non-negativity inequalities $\Sigma_T\geq 0$, we have that, for $s>m+1$, only the vector of all zero is possible as $d$, i.e. the cone collapses to zero and we obtain (i).

Now, summing the $(m,s)$-inequalities over all $i\not= k$, we obtain $(s+1)(\Sigma_T-d_{x_k})\leq (m+1)\Sigma_T$, i.e. $(s+1)d_{x_k}\geq (s-m)\Sigma_T$. This inequality and the $(m,s)$-simplex equality gives
\begin{equation*}
(s-m)\Sigma_T \leq (s+1)d_{x_k}\leq \Sigma_T.\mbox{~~~~~~~~~~~~~~}(*)
\end{equation*}
For $s=m+1$, this inequality imply $(m+2)d_{x_k}=\Sigma_T$, i.e. $d$ is a positive multiple of the vector of all ones and we obtain (ii).

Let $m\leq s\leq m+1$. Then, the inequalities $(*)$ imply the inequality $\Sigma_T\geq 0$ and, therefore, the non-negativity inequalities. This gives (iii).
\\

We will now assume $0<s<m+1$.

\begin{remark}
We got by computations the following new facts about small cones $MET_n$:

(i) the full symmetry group of $MET_n$ is $Sym(n)$ for $n=3, 5, \dots, 14$, and for $n=4$ it is $Sym(4)\times Sym(3)$;

(ii) the diameter of $G(MET_7)$ is three;

(iii) using the list of $1550825600$ vertices of the metric polytope $MET^{\Box}_8$ (see \cite{DFPS}); we obtained $3918$ orbits of extreme rays of the metric cone $MET_8$;

(iv) the cone $MET_7$ has $46$ orbits of extreme rays and not $41$, as was given in \cite{GrMet} and \cite{DL}.
\end{remark}

Note, that the group $Sym(n)$ acts in an obvious way on $V_n$ and so, on $SMET^{m,s}_n$, which proves that $Sym(n)$ is a subgroup of the full symmetry group of $SMET^{m,s}_n$. In order to find the symmetry group of $SMET^{m,s}_n$,
 one should consider the group of automorphisms of the graphs $G(SMET^{m,s}_n)$ and $G({SMET^{m,s}_n}^*)$. If one of those groups is equal to $Sym(n)$, then the full automorphism group of the cone $SMET^{m,s}_n$ is $Sym(n)$

%When those groups can be computed (for example using the {\it nauty} program\cite{MK}) and the 

We checked by computer, using the {\it nauty} program (\cite{MK}), for $(m,s,n)=(2,2,5)$, $(2,1,1)$, $(2,1,6)$, $(3,2,5)$, $(3,3,6)$, $(3,1,6)$, $(3,2,6)$, $(4,1,7)$ and for $MET_n$ with $5\leq n \leq 14$, the validity of the following conjecture:
\begin{conj}
The symmetry group of $SMET^{m,s}_n$ is $Sym(n)$, the only exception being $MET_4=SMET^{1, 1}_4$, for which it is $Sym(4)\times Sym(3)$.
\end{conj}

Remark that in the $SCUT^{m,s}_n$ case, the situation is more involved; for example the symmetry group of $SCUT^{3,3}_6$ has size $518400$ and not $720$ as we expected.

%$SCUT^{2, 2}_5$, $HCUT^2_5$, $SCUT^3_6$, $HCUT^{4}_7$, $HCUT^2_6$ 

Remark that 

%The unsolved cases are $SCUT^{3, 3}_6$, $SMET^{4,3}_7$, $SMET^{4,2}_7$.

\begin{theor}
For the facets of $SMET^{m,s}_n$ holds:

(i) The non-negativity facet $NN_{A-\{i\}}$ is non-adjacent to the $(m,s)$-simplex facet $ST_{A; i}$ in the cone $SMET^{m, s}_n$,

(ii) two $m$-simplex facets are non-adjacent if they have the same support,

(iii) for  $m-1\leq s<m$, the non-negativity facets $NN_{A}$ and $NN_{B}$ are non-adjacent if $|A\cap B|=m$.
\end{theor}
\proof (i) If $d$ is a $m$-hemi-metric satisfying $d(A-\{i\})=0$ and $sd(A-\{i\})=\sum_{k\in A-\{i\}} d(A - \{k\})$, then we have $d(A-\{k\})=0$ for all $k\in A-\{i\}$ and $d$ lies in a space of insufficient rank.

(ii) Assume that $d$ is incident to the facets $ST_{T, i}$ and $ST_{T, j}$ with $T$ being a subset of $V_n$ of size $m+2$ and $i,j\in T$. Then we have $2d_{i}=\Sigma_T$ and $2d_{j}=\Sigma_T$, which implies $d_i=d_j$ and $0=\sum_{k\in T-\{i,j\}} d(T-\{k\})$. So, $d(A-\{k\})=0$ and the rank is again too low.

(iii) If $|A\cap B|=m$, then write $A\cup B=(A\cap B)\cup \{i,j\}$ and set
$d_{(A\cup B)-\{k\}}=d_k$ for $k\in A$.
The non-negativity inequalities $NN_A$ and $NN_B$ have the form $d_{i}\geq 0$ and $d_{j}\geq 0$. Let $d$ lie on $NN_{A}$ and $NN_B$. Then $d_{i}=d_{j}=0$. Summing $(m,s)$-simplex inequalities $(s+1)d_k\leq \Sigma_{A\cup B}$ over all $k\not= i,j$, and taking in attention that $d_{i}=d_{j}=0$, we obtain $(s+1)\Sigma_{A\cup B}\leq m\Sigma_{A\cup B}$. If $m-1<s<m$, this inequality holds only if $\Sigma_{A\cup B}=0$. Hence $d_{k}=0$ for all $k\in A\cup B$, i.e. the rank is too low.

If $s=m-1$, the last inequality is equality. This means that all summed inequalities are equalities, too. Hence $md_{k}=\Sigma_{A\cup B}$ for all $k\not= i,j$. This implies, that the intersection of $NN_A$ and $NN_B$ is the ray $d_i=d_j=0$, $d_k=const\geq 0$, $k\not= i,j$. Hence codimension of this intersection is $m+1$, which is strictly greater than two, needed for adjacency.

\begin{conj}
For the facets of $SMET^{m,s}_n$ holds:

(i) two $(m,s)$-simplex facets are not adjacent if and only if $s=1$ and they have the same support,

(ii) a non-negativity and a $(m,s)$-simplex facets are not adjacent if and only if they are {\em conflicting} (i.e. there exist a position, for which they have non-zero values of different sign),

(iii) two non-negativity facets, say, $NN_A$ and $NN_B$ are not adjacent if and only if $|A\cap B|=m$ and $m-1\leq s<m$,

(iv) the ridge graph is complete if $s\geq m$; otherwise, it has diameter two and, for $1<s<m-1$, it is $K_{(n-m){n \choose m+1}}-{n\choose m+1} K_2$.
\end{conj}
This conjecture was checked for all cases of Table \ref{tab:MainLovelyTable}.

%$(m,s,n)=(2,1, 5)$, $(2,3/2, 5)$, $(2, 2, 5)$, $(2, 1, 6)$, $(2, 1, 7)$, $(3, 2, 6)$, $(3, 3/2,6)$, $(3, 3, 6)$, $(3, 1, 6)$, $(3, 1, 7$, $(4, 4, 7)$, $(4, 3, 7)$, $(4, 1, 7)$, $(5, 1, 7)$, $(5, 2, 8)$, $(5, 3, 8)$, $(5, 4, 8)$, $(5, 5, 8)$.
%These conjectures were checked for $SMET^{2, x}_5$ ($x=1, 1.5, 2$), $HMET^2_6$, $SMET^{3, 2}_6$, $SMET^{3, 3/2}_6$, $SMET^{3, 3}_6$, $HMET^3_6$, $SMET^{4, 4}_7$, $SMET^{4, 3}_7$, $HMET^4_7$, $HMET^5_7$, $SMET^{5, 2}_8$, $SMET^{5, 3}_8$, $SMET^{5, 4}_8$, $SMET^{5,5}_8

\begin{proposition}
The number of orbits of $0,1$ valued extreme rays and the minimal number of zeros of an extreme ray are:
\begin{center}
\small
\begin{tabular}{|c|c|c|c|}
\hline
$(3,5)$ for $HMET^{2}_5$&$(2,0)$ for $SMET^{2, 2}_5$&&\\
\hline
$(4,8)$ for $HMET^3_6$&$(5,3)$ for $SMET^{3,2}_6$ &$(2,0)$ for $SMET^{3,3}_6$&\\
\hline
$(5,13)$ for $HMET^{4}_7$&$(10,7)$ for $SMET^{4, 2}_7$&$(7,3)$ for $SMET^{4, 3}_7$&$(2,0)$ for $SMET^{4, 4}_7$\\
\hline
$(6,9)$ for $HMET^{2}_6$&$(3,0)$ for $SMET^{2, 2}_6$&&\\
\hline
\end{tabular}
\end{center}
\end{proposition}

%\newline
By analogy with the binomial coefficient ${n \choose m}$, let us denote by ${A \choose m}$ the {\em set of all $m$-subsets of the set $A$}.

\begin{definition}
(i) For any function $d$ on the set ${V_n \choose m+1}$ one can define a function $d^{ze}$ on the set ${V_{n+1} \choose m+2}$, called {\it zero-extension}, by
\begin{equation*}
d^{ze}_{S}=\left\lbrace\begin{array}{cl}
0&\mbox{~if~}S\subset V_{n}\\
d_{S\cap V_{n}},&\mbox{~otherwise,}
\end{array}\right.
\end{equation*}
for any $(m+2)$-subset $S$ of $V_{n+1}$.

(ii) For any function $d$ on the set ${V_n \choose m+1}$ one can define a function $d^{vs}$ on the set ${V_{n+1} \choose m+1}$, called {\em vertex-splitting}, (of the vertex $n$ into two vertices $n$ and $n+1$) by
\begin{equation*}
d^{vs}_{S}=\left\lbrace\begin{array}{cl}
0&\mbox{~if~}\{n,n+1\}\subset S,\\
d_S&\mbox{~if~}S\subset V_n,\\
d_{(S-\{n+1\})\cup \{n\}}&\mbox{~if~}n+1\in S\mbox{~and~}n\notin S,
\end{array}\right.
\end{equation*}
for any $(m+1)$-subset $S$ of $V_{n+1}$.
\end{definition}

For example, for the cuts $\delta_{\{1\}, \{2, 3, 4\}}$ and $\delta_{\{1, 2\}, \{3, 4\}}$ (which are representatives of two orbits of extreme rays of $MET_{4}$) zero-extension are vectors 
\begin{center}
\scriptsize
\begin{tabular}{cccccccccc}
 $\overline{45}$&$\overline{35}$&$\overline{34}$&$\overline{25}$
&$\overline{24}$&$\overline{23}$&$\overline{15}$&$\overline{14}$
&$\overline{13}$&$\overline{12}$\\
\hline
0&0&1&0&1&1&0&0&0&0\\
0&0&0&0&1&1&0&1&1&0\\
\hline
\end{tabular}
\end{center}
(which are extreme rays of $HMET^{2}_5$) and vertex-splitting are vectors
\begin{center}
\scriptsize
\begin{tabular}{cccccccccc}
12&13&14&15&23&24&25&34&35&45\\
\hline
1&1&1&1&0&0&0&0&0&0\\
0&1&1&1&1&1&1&0&0&0\\
\hline
\end{tabular}
\end{center}
i.e. $\delta_{\{1\}, \{2, 3, 4, 5\}}$ and $\delta_{\{1, 2\}, \{3, 4, 5\}}$ (which are extreme rays of $MET_{5}$). 

The graph $G_{d^{ze}}$ is equal to $G_{d}$ for any $0,1$ valued vector. The graph $G_{d^{vs}}$ is obtained from $G_{d}$ by splitting each 
vertex, associated to a set $S\subset V_n$, with $n\in S$, in two vertices: one for subset $S$ and another one for the subset $(S-\{n\})\cup \{n+1\}$.

\begin{theor}\label{zero-extension-vertex-splitting}
(i) Zero-extension of any extreme ray of $SMET^{m,s}_{n}$ is an extreme ray of $SMET^{m+1,s}_{n+1}$;

(ii) vertex-splitting of any extreme ray of $HMET^{m}_n$ is an extreme ray of $HMET^{m}_{n+1}$.
\end{theor}
%Part (i) was checked for $(m,s,n)=(4,1,7)$, $(4,2,7)$, $(4,3,7)$, $(4,4,7)$, $(2,2,5)$, $(2,2,6)$ and $(3,3,6)$. Part (ii) is true for $m=1$ (see \cite{GrMet}) and was checked for $HMET^{2}_5$,  $HMET^{3}_6$, $HMET^{4}_7$, $HMET^{5}_8$,  $HMET^{2}_6$.
\proof (i) If $d$ is an extreme ray of $SMET^{m,s}_n$, then one can check easily the validity of $(m+1,s)$-simplex inequalities and non-negativity inequalities for $d^{ze}$.

The ray $d^{ze}$ is incident to the non-negativity facet $NN_A$ if $A\subset V_n$ and to the $(m,s)$-simplex facet $ST_{A;i}$, if $n+1\in A$ and $d$ is incident to $ST_{A-\{n+1\};i}$.

Assume now that $e$ is a ray of $SMET^{m+1,s}_{n+1}$, which is incident to all facets incident to $d^{ze}$. Then we obtain $e(A)=0$ if $A\subset V_n$;

so, the restriction of $e$ to the subsets, containing $n+1$, is identified to a function on ${V_n \choose m+1}$, which is incident to the facets incident to $d$. So, we get $e=\lambda d^{ze}$ with $\lambda\geq 0$ by non-negativity.

(ii) If $d$ is an extreme ray of $HMET^{m,s}_n$, then one can check easily the validity of $(m,s)$-simplex inequalities and non-negativity inequalities for $d^{vs}$.

The ray $d^{vs}$ is incident to the non-negativity facet $NN_A$ if $\{n,n+1\}\subset A$ and to the $(m,s)$-simplex facet $ST_{A;i}$, if $A\subset V_n$ and $d$ is incident to $ST_{A;i}$.

Assume now that $e$ is a ray of $HMET^{m,s}_{n+1}$, which is incident to all facets incident to $d^{vs}$. Then one obtain $e(A)=0$ if $\{n,n+1\}\subset A$.

Now, if $\{n,n+1\}\subset T$, then, applying $m$-simplex-inequalities $ST_{T;n+1}$ and $ST_{T; n}$, we get $e(T-\{n\})=e(T-\{n+1\})$.

Since the restriction of $e$ on $V_n$ yield a vector, which is incident 
to all facets incident to $d$, one obtain (since $d$ is an extreme ray)
that the restriction $e_{|V_n}$ is a multiple of $d$. So, we get
$e_{|V_n}=\lambda d$ with $\lambda \geq 0$. Above equalities yield
$e=\lambda d^{vs}$.\\

The same result holds for the similar notion of vertex-splitting of an ray in $QMET_n$ (see \cite{DDP}).

\begin{theor}\label{Theorem-completely-solvable}
For the cone $SMET^{m,s}_{m+2}$ holds:

(i) It has only one orbit of extreme rays. Each extreme ray contains a vector with $\lfloor s\rfloor+1$ components $1$, one component $s-\lfloor s\rfloor$ and the other ones $0$; all such vectors appear on different extreme rays.

(ii.1) If $s$ is integer, the skeleton is the Johnson graph $J(m+2,s+1)$;

(ii.2) if $s$ is not integer, then two extreme rays are adjacent if and only if they either have the same support, or they differ only by the position of the value $s-\lfloor s\rfloor$ in the associated vector.

(iii.1) If $m\leq s<m+1$, then both, the skeleton and the ridge graph, are $K_{m+2}$;

(iii.2) if $1< s < m-1$, then the ridge graph is $K_{(m+2)\times 2}$;

(iii.3) if $s=1<m-1$ or $1<s=m-1$, then it is $K_{(m+2)\times 2}-K_{m+2}$;

(iii.4) if $s=1=m-1$, then it is $K_{(2+2)\times 2}-2K_{2+2}$, (i.e. $3-\mbox{cube}$).

\end{theor}
\proof (i) The cone is defined by $\sum d_i\geq (s+1)d_k\geq 0$ for all $1\leq k\leq m+2$.

In fact, $SMET^{m, s}_{m+2}$ has dimension $m+2$ and two orbits of facets, each consisting of $m+2$ linearly independent members: $(m,s)$-simplex facets - all $(1, -s)$-vectors of length $m+2$ with exactly one $-s$, and non-negativity facets - all $(0, 1)$-vectors of length $m+2$ with exactly one $1$.
Fix an extreme ray $d$ of the cone. It lies on $m+1$ linearly independent facets. Without loss of generality, one can suppose that it lies on $(m,s)$-simplex facets with $-s$ on positions $1, \dots, p$ only and on non-negativity facets with $1$ on positions $p+2, \dots, m+2$ only. So, $d_1=\dots=d_p=t$, say; $d_{p+2}=\dots=d_{m+2}=0$ and $(-s+p-1)t+d_{p+1}=0$. The validity of the $(m,s)$-simplex facet with $-s$ on position $p+1$ and of the non-negativity facets imply $t> 0$, $d_{p+1}\geq 0$, and $pt-sd_{p+1}\geq 0$. This yield the inequalities $p-1\leq s\leq p$; if $s$ is integer then the values $p=s+1$ and $p=s$ yield the same solution; so, one can assume, in general, $p=\lfloor s\rfloor +1$ and we are done.

(ii) If $s$ is integer, the vertices of the skeleton are the same as for the Johnson graph $J(m+2, s+1)$. Let us see now that it is this Johnson graph, i.e. two extreme rays are adjacent if and only if the rank of the set of facets, containing them both, is $m$. If corresponding vectors have, say, $i$ common ones, then this rank is $m+2i-2s$ (namely, $i$ simplex-facets and $m+i-2s$ non-negativity facets). The maximum of this number is $m$ and is attained exactly for $i=s$. The diameter of skeleton is $\min(s+1, m-s+1)$.

If $s$ is non-integer, then each extreme ray belongs to exactly $m+1$ (linearly independent) facets: $\lfloor s\rfloor+1$ $(m,s)$-simplex facets with $-s$ on position, where the ray has $1$, $m-\lfloor s\rfloor$ non-negativity facets with $1$ on position, where the ray has $0$. So, the adjacency of extreme rays follows.

(iii) First, non-negativity inequalities are not facets if and only if 
$m\leq s<m+1$. Two $(m,s)$-simplex facets are adjacent, unless $s=1<m$. 
Two non-negativity facets are adjacent, unless $s\geq m-1$. Fixed 
$(m,s)$-simplex and non-negativity facets are adjacent if and only if 
there is no position, in which the first has $-s$ and the second has $1$. 
So, the ridge graph follows.

%(In the last case, those facets are, moreover, partition the set of extreme 
%rays into two parts of equal size). 

\begin{remark}
(i) Any $SMET^{1,s}_{3}$ has only one orbit of extreme rays, represented by
$(0,s,1)$ for $0<s\leq 1$ and by $(s-1,1,1)$ for $1\leq s\leq m+1=2$. Both, 
the skeleton and the ridge graph of this cone, are $C_6$ for $0<s<1$ and $K_3$ for $1\leq s<2$.

(ii) The number of extreme rays (number of orbits) of $SMET^{1,s}_n$ is
$54(5)$, $2900(35)$, $988105(1567)$ for $s=\frac{1}{2}$ and $n=4,5,6$; it is
$25(4)$, $1235(24)$, $530143(890)$ for $s=\frac{3}{2}$ and $n=4,5,6$.
\end{remark}

%The case with $s$ non-integer is also solvable:
%\begin{conj}
%The cone $SMET^{m,s}_{m+2}$ with $s=k+\delta$, $k\in\{1, 2, \dots, m-1\}$, $0<\delta<1$ has one orbit of extreme rays. A vector lying on one of these extreme ray has $k+2$ non-negative components, $k+1$ of its component are $1$ and the last is $\delta$.\\
%There are two cases of adjacency between those extreme rays:
%\begin{enumerate}
%\item They have the same support;
%\item they differ by the position of the $\delta$ in the associated vector.
%\end{enumerate}
%\end{conj}
We will suppose from now that $n\geq m+3$.

The {\it incidence number} of a facet (or of an extreme ray) is the number of
extreme rays lying on this facet (or, respectively, of facets containing
this extreme ray). The {\it adjacency number} of a facet or extreme ray is the
number of adjacent facets or extreme rays. In Tables below, we
give representative of orbits of facets (or extreme rays) with adjacency, 
incidence of their representatives and orbit-size.

The {\it representation matrix} of skeleton (or ridge) graph is the 
square matrix, where on the place $i,j$ we put the number of members of 
orbit $O_j$ of extreme rays (or facets, respectively), which are adjacent 
to a fixed representative of orbit $O_i$.

\section{$m$-partitions hemi-metrics}

Given an unordered $(m+1)$-partition $S_{1}, \ldots, S_{m+1}$ of 
$V_n = \{1, 2, \ldots, n\}$, the $m$-hemi-metric, which we denote
$\alpha(S_{1}, \ldots, S_{m+1})$, and call {\it $m$-partition hemi-metric} is
defined by 
\begin{equation*}
\alpha(S_{1}, \ldots, S_{m+1}) (x_1, \dots, x_{m+1})=
\left\lbrace\begin{array}{cl}
1&\mbox{~if~}x_1 \in S_1, \dots, x_{m+1} \in S_{m+1}\\
0&\mbox{~otherwise.}
\end{array}\right.
\end{equation*}
It is easy to see that $\alpha(S_{1}, \ldots, S_{m+1})$ is a $m$-hemi-metric, and, for $m=1$, it is the usual cut semi-metric.

The graph $G_v=K_{|S_1|}\times K_{|S_2|}\times\dots\times K_{|S_{m+1}|}$ is associated to the $m$-partition hemi-metric vector $v=\alpha(S_1, S_2, \dots, S_{m+1})$; the zero-extension of $v$ is the $(m+1)$-partition hemi-metric $\alpha(S_1, \dots, S_{m+1},\{n+1\})$, while its vertex-splitting is the $m$-partition hemi-metric $\alpha(S_1, \dots, S_p\cup \{n+1\}, \dots,  S_{m+1})$ if $n\in S_p$.

\begin{theor}
The $m$-partition hemi-metrics are extreme rays of $HMET^m_n$.
\end{theor}
\proof Any ray $\alpha(S_1, S_2, \dots, S_{m+1})$ is a function on ${V_n\choose m+1}$. Using symmetry and above remark, it can be viewed as a vertex-splitting of ray $\alpha(\{1\}, \{2\}, \dots, \{m+1\})$ on ${V_{m+1}\choose m+1}$. The hemi-metric cone $HMET^m_{m+1}$ has dimension one and $\alpha(\{1\}, \{2\}, \dots ,\{m+1\})$ generates it; so, applying Theorem \ref{zero-extension-vertex-splitting}, we obtain the result.

%This conjecture holds for $n=m+2$ (it is implied by theorem \ref{Theorem-completely-solvable}) and we checked it for $(m, n)=(2, 5)$, $(3, 6)$, $(4, 7)$, $(5, 8)$, $(6, 9)$, $(2, 6)$, $(3, 7)$, $(4, 8)$, $(2, 7)$, $(3, 8)$.\\

\begin{definition}
Denote by $HCUT_{n}^m$ the cone generated by all $m$-partition hemi-metrics on $V_n$.
\end{definition}

For $m=1$ this cone is $CUT_n$. Clearly, $HCUT_n^m \subseteq HMET_n^m$ with equality for $n=m+2$.

\begin{conj}
Let $m> 1$ and $n\geq m+2$; then $HCUT^m_n=HMET^m_n$ only if $n=m+2$.
\end{conj}

\begin{conj}\label{Conjecture-adjacencies-extrays}
In the skeleton of $HCUT^m_n$,  two $m$-partition hemi-metrics are 
non-adjacent if and only if, up to permutations, corresponding 
$(m+1)$-partitions  can be written as $(S_1, \dots, S_{m+1})$ and 
$(S'_1, \dots, S'_{m+1})$, which differ only by $S_1, S_2, S_3=A\cup B, C, D$
 and $S'_1, S'_2, S'_3=A, B, C\cup D$ for some disjoints sets $A, B, C, D$.
%Two extreme rays of $HCUT^m_n$ are non-adjacent iff, up to permutations,
%corresponding (m+1)-partitions of $n$-set can be written as $(S_1, \dots, S_{m+1})$ and $(S'_1, \dots, S'_{m+1})$, which differ only bt $S_1=S'_1\cup S'_2$ and $S_2\cup S_3=S'_3$.
%With this criterion, the orbit of extreme ray, having all $S_i$ of size one, except, say $S_1$, form a dominating clique and so the diameter is $2$ or $3$.
\end{conj}
This conjecture holds for $n=m+2$ and it was checked also for $(m,n)=(2, 5)$, $(3, 6)$, $(4, 7)$, $(5, 8)$, $(6, 9)$, $(2, 6)$, $(3, 7)$, $(4, 8)$, $(2, 7)$.

If this conjecture holds, it will imply, that the skeleton graph of $HCUT^m_n$ has diameter two, since any two non-adjacent $m$-partitions hemi-metrics are both adjacent to $m$-partition hemi-metric $\alpha(A\cup C, B, D, S_4, \dots, S_{m+1})$.

\begin{conj}
(i) Adjacency rule in $HMET^m_n$ is the same for the first orbit, represented by $\{1\}, \{2\}, \dots, \{m\}, \{m+1, \dots, n\}$,

(ii) but two elements of the second orbit (represented by $\{1\}, \{2\}, \dots, \{m-1\}, \{m, m+1\}, \{m+2, \dots, n\}$) are non-adjacent if and only if the rule as above is completed by the condition 

$rank(\inf(a,b))-(m+1)=(m+1)-rank(\sup(a,b))>1$ (rank here is the number of parts, 
while $\inf(a,b)$ and $\sup(a,b)$ are join and union operation in the lattice of partitions).
\end{conj}

\begin{conj}
The non-negativity inequalities and the $m$-simplex inequalities are facets of $HCUT^m_n$.
\end{conj}
This conjecture was verified for $(m,n)=(3, 6)$, $(4, 7)$, $(5, 8)$ and for $(2,6)$.

\subsection{Facets of $HCUT^{2}_6$, $HCUT^{4}_7$}
\noindent We list in Table \ref{InfoConcerningTheConeHCUT26} and \ref{InfoConcerningTheConeHCUT47} only the orbits with incidence greater or equal to $29$, besides the non-negativity and $m$-simplex facets.

\begin{table}
\begin{center}
\tiny
\begin{tabular}{|p{5pt}p{5pt}p{4pt}p{5pt}p{4pt}p{4pt}p{5pt}p{4pt}p{4pt}p{4pt}p{4pt}p{5pt}p{4pt}p{5pt}p{4pt}p{4pt}p{5pt}p{5pt}p{5pt}p{7pt}|c|c|c|}
\hline
123&124&125&126&134&135&136&145&146&156&234&235&236&245&246&256&345&346&356&456&Adj.&Size&Inc.\\
\hline
$-$1&0&0&1&0&1&0&0&0&1&0&1&2&0&0&$-1$&0&0&0&0&15270&360&45\\
$-1$&$-1$&0&2&2&0&1&0&1&0&2&0&1&0&1&0&0&$-2$&0&0&5908&180&42\\
$-1$&$-1$&$-1$&3&$-1$&1&1&1&1&1&$-1$&1&1&1&1&1&1&1&$-1$&$-1$&600&180&33\\
$-1$&$-1$&$-1$&3&1&1&1&1&1&1&1&1&1&1&1&1&$-1$&$-1$&$-1$&$-1$&939&60&33\\
$-1$&$-1$&1&1&0&0&1&1&0&0&0&0&1&2&1&$-1$&1&$-1$&1&0&1579&720&33\\
$-1$&$-1$&0&2&1&$-1$&1&1&1&0&1&$-1$&1&1&1&0&3&$-1$&1&$-1$&496&360&32\\
$-1$&$-1$&1&1&2&0&1&1&0&0&2&1&0&0&1&0&$-1$&$-1$&0&0&1856&180&32\\
$-1$&$-1$&1&1&0&0&1&0&1&1&1&1&1&1&1&$-1$&$-1$&0&0&0&515&360&31\\
$-1$&0&0&1&0&0&1&1&$-1$&1&0&0&1&1&1&$-1$&2&0&0&0&629&360&31\\
$-1$&$-1$&0&2&1&$-1$&1&1&1&0&2&$-1$&2&0&1&1&1&$-2$&1&0&404&720&30\\
$-1$&$-1$&1&1&0&0&1&1&0&0&0&1&0&2&1&0&0&0&1&$-1$&458&720&30\\
$-1$&$-1$&1&1&2&0&1&0&1&1&3&1&1&1&1&$-1$&$-1$&$-2$&0&0&558&360&30\\
$-1$&$-1$&2&2&2&$-1$&2&2&$-1$&$-1$&2&2&$-1$&$-1$&2&$-1$&$-1$&$-1$&2&2&2265&12&30\\
$-1$&0&0&1&0&0&1&1&$-1$&1&0&1&2&0&0&$-1$&1&1&0&0&441&720&30\\
$-1$&0&0&1&0&1&0&1&1&0&1&0&0&1&0&1&0&1&1&$-2$&1867&120&30\\
$-1$&$-1$&$-2$&4&$-2$&2&1&2&1&2&$-2$&2&1&2&1&2&4&2&$-2$&$-2$&149&180&29\\
$-1$&$-1$&0&2&1&$-1$&1&1&1&0&1&$-1$&1&1&1&0&2&$-2$&2&0&473&360&29\\
$-1$&$-1$&0&2&1&2&0&2&0&2&2&1&2&1&2&$-2$&$-1$&$-2$&0&0&323&360&29\\
$-1$&$-1$&0&2&3&$-1$&1&0&2&1&3&1&1&$-1$&1&0&1&$-3$&1&0&288&720&29\\
$-1$&$-1$&1&1&$-1$&1&1&1&1&$-1$&0&0&1&0&1&1&1&0&0&0&283&360&29\\
\hline
\end{tabular}
\end{center}
\caption{Representatives of orbits of facets of $HCUT^{2}_6$ with incidence at least $29$}
\label{InfoConcerningTheConeHCUT26}
\end{table}

\begin{table}
\begin{center}
\tiny
\begin{tabular}{|p{5pt}p{5pt}p{5pt}p{5pt}p{5pt}p{5pt}p{5pt}p{5pt}p{5pt}p{5pt}p{5pt}p{5pt}p{5pt}p{5pt}p{5pt}p{5pt}p{5pt}p{5pt}p{5pt}p{5pt}p{5pt}|c|c|c|}
\hline
&&&&&&&&&&&&&&&&&&&&&&&\\[-1mm]
$\overline{67}$&$\overline{57}$&$\overline{56}$&$\overline{47}$&$\overline{46}$&$\overline{45}$&$\overline{37}$&$\overline{36}$&$\overline{35}$&$\overline{34}$&$\overline{27}$&$\overline{26}$&$\overline{25}$&$\overline{24}$&$\overline{23}$&$\overline{17}$&$\overline{16}$&$\overline{15}$&$\overline{14}$&$\overline{13}$&$\overline{12}$&Adj.&Size&Inc.\\
\hline
$-1$&0&1&0&1&0&0&1&0&0&1&0&1&1&1&2&1&0&0&0&$-1$&3490&420&50\\
$-1$&0&1&0&1&0&0&1&0&0&1&0&1&1&1&1&2&$-1$&1&1&0&343&1260&37\\
$-1$&$-1$&2&0&1&1&0&1&1&0&2&1&1&0&0&2&1&1&2&2&$-2$&320&1260&36\\
$-1$&$-1$&2&1&0&2&1&0&2&0&1&2&0&0&0&1&2&0&0&0&0&353&630&34\\
$-1$&0&1&0&1&0&1&0&1&1&1&0&1&1&0&1&2&$-1$&1&0&0&429&1260&34\\
$-1$&0&1&0&1&0&1&0&1&1&1&2&$-1$&1&0&1&2&1&$-1$&0&2&405&840&34\\
$-1$&$-1$&2&$-1$&2&2&0&1&1&1&2&1&1&1&0&2&1&1&1&2&$-2$&66&840&32\\
$-1$&$-1$&2&0&1&1&1&0&2&1&1&2&0&1&2&2&1&1&0&$-1$&$-1$&213&1260&32\\
$-1$&0&1&0&1&0&1&0&1&1&1&2&$-1$&1&0&2&1&2&0&$-1$&1&394&2520&32\\
$-1$&$-1$&2&$-1$&2&2&$-1$&2&2&2&2&1&1&1&1&2&1&1&1&1&$-2$&30&105&30\\
$-1$&$-1$&2&$-1$&2&2&1&1&1&1&1&1&1&1&$-1$&1&1&1&1&$-1$&2&33&420&29\\
$-1$&$-1$&2&$-2$&3&3&1&2&2&1&2&1&1&2&$-1$&2&1&1&2&3&$-2$&44&2520&29\\
$-1$&$-2$&3&0&1&2&1&2&1&1&2&1&2&0&$-1$&2&1&2&2&3&$-2$&89&2520&29\\
\hline
\end{tabular}
\end{center}
\caption{Representatives of orbits of facets of $HCUT^{4}_7$ with incidence at least $29$}
\label{InfoConcerningTheConeHCUT47}
\end{table}

\subsection{The cone $HMET^2_6$}
\begin{table}
\begin{center}
\tiny
\begin{tabular}{|c|p{4pt}p{4pt}p{4pt}p{4pt}p{4pt}p{4pt}p{4pt}p{4pt}p{4pt}p{4pt}p{4pt}p{4pt}p{4pt}p{4pt}p{4pt}p{4pt}p{4pt}p{4pt}p{4pt}p{6pt}|c|c|c|}
\hline
&123&124&125&126&134&135&136&145&146&156&234&235&236&245&246&256&345&346&356&456&Adj.&Size&Inc.\\
\hline
F1&$-1$&0&0&1&0&0&1&0&0&0&0&0&1&0&0&0&0&0&0&0&75&60&4001\\
F2&0&0&0&0&0&0&0&0&0&0&0&0&0&0&0&0&0&0&0&1&67&20&3939\\
\hline
E1&0&0&0&0&0&0&0&0&0&1&0&0&0&0&0&1&0&0&1&1&2778&15&64\\
E2&0&0&0&0&0&0&0&0&1&1&0&0&0&0&1&1&0&1&1&0&1321&60&56\\
E3&0&0&1&1&1&0&1&1&0&0&1&1&0&0&1&0&0&0&1&1&1030&12&40\\
E4&0&0&0&0&0&1&1&1&1&0&0&1&1&1&1&0&0&0&0&0&818&15&48\\
E5&0&0&0&0&0&0&1&1&0&1&0&0&1&1&0&1&1&1&0&0&731&180&48\\
E6&0&0&0&1&0&1&0&1&1&0&0&1&1&1&0&1&0&1&1&0&358&180&40\\
E7&0&0&1&1&1&0&1&1&0&0&1&1&0&0&1&0&0&2&1&1&270&120&36\\
E8&0&0&1&1&1&0&1&1&0&2&1&1&2&2&1&0&0&0&1&1&93&120&28\\
E9&0&0&1&1&1&0&1&1&0&2&1&1&0&2&1&0&0&2&1&1&66&240&28\\
E10&0&0&1&1&1&0&1&1&0&2&1&1&0&0&1&2&0&2&1&1&51&360&28\\
E11&0&0&1&1&1&0&1&0&1&1&1&1&0&1&0&1&1&1&0&2&47&120&28\\
E12&0&0&0&1&0&0&1&0&1&2&0&0&2&0&2&1&0&2&1&1&46&60&39\\
E13&0&0&0&1&0&0&1&1&1&2&0&0&2&1&2&1&1&2&1&0&37&360&31\\
E14&0&0&0&1&0&0&1&2&2&2&0&0&2&2&1&1&2&1&1&0&37&180&31\\
E15&0&0&0&2&0&0&2&2&1&1&0&0&2&2&1&1&2&1&1&0&37&60&31\\
E16&0&0&0&2&0&3&1&3&1&1&0&3&1&3&1&1&0&2&1&1&32&90&28\\
E17&0&0&0&1&0&1&1&2&2&2&0&1&2&2&1&1&1&1&0&0&30&720&27\\
E18&0&0&0&2&0&1&2&2&1&1&0&1&2&2&1&1&1&1&0&0&30&360&27\\
E19&0&0&0&1&0&1&1&1&1&2&0&1&2&1&2&1&2&2&0&0&30&360&27\\
E20&0&0&0&2&0&2&1&2&1&1&0&2&1&2&1&1&2&2&0&0&30&180&26\\
E21&0&0&0&2&0&1&1&1&1&2&0&1&1&1&1&2&2&2&0&0&30&180&27\\
E22&0&0&0&2&0&2&1&3&1&1&0&2&1&3&1&1&1&2&0&1&29&360&26\\
E23&0&0&0&2&0&1&1&2&1&1&0&1&1&2&1&1&1&2&3&0&29&360&26\\
E24&0&0&1&1&1&0&1&0&2&2&1&1&2&1&1&0&1&0&1&2&27&360&23\\
E25&0&0&0&1&1&1&1&2&2&2&1&1&2&2&1&1&0&0&0&0&27&360&27\\
E26&0&0&1&1&2&0&1&0&1&2&2&1&2&1&2&0&2&0&1&1&27&360&23\\
E27&0&0&1&2&3&0&1&3&1&1&3&1&1&2&1&0&0&1&2&1&27&360&23\\
E28&0&0&1&2&2&1&2&1&2&0&2&2&0&2&0&1&0&2&1&1&27&180&23\\
E29&0&0&0&2&1&1&2&2&1&1&1&1&2&2&1&1&0&0&0&0&27&180&27\\
E30&0&0&0&2&2&2&1&2&1&1&2&2&1&2&1&1&0&0&0&0&27&60&26\\
E31&0&0&1&2&1&1&1&2&1&0&1&2&1&1&1&1&0&3&0&3&26&720&23\\
E32&0&0&0&2&1&1&1&2&1&1&1&1&1&2&1&1&0&3&3&0&26&180&24\\
E33&0&0&1&1&2&0&1&1&1&2&2&1&2&2&2&0&1&0&1&0&25&720&23\\
E34&0&0&1&1&1&0&2&1&1&2&1&1&1&2&2&0&2&0&2&0&25&720&23\\
E35&0&0&1&1&0&1&2&2&1&0&0&2&1&1&2&0&1&3&1&1&25&360&25\\
E36&0&0&1&1&0&1&2&2&1&0&0&2&1&1&2&0&1&1&1&3&25&360&25\\
E37&0&0&1&1&1&0&1&2&2&2&1&1&2&1&1&0&1&0&1&0&25&360&23\\
E38&0&0&1&1&1&0&2&2&1&1&1&1&1&1&2&3&1&0&1&0&23&720&22\\
E39&0&0&1&2&1&1&1&2&0&1&1&2&1&1&2&0&0&2&3&1&22&720&21\\
E40&0&0&1&2&1&1&2&2&1&0&1&2&0&1&3&1&0&2&1&1&22&720&21\\
E41&0&0&1&2&2&1&1&1&1&1&2&2&1&2&1&0&0&0&3&3&21&360&21\\
\hline
\end{tabular}
\end{center}
\caption{Representatives of orbits of facets and extreme rays of $HMET^{2}_6$}
\end{table}
The graph $G_v$ for orbits $E_1$, $E_2$, $E_3$, $E_4$, $E_5$ and $E_6$ is $K_4$, $\overline{C_6}$ (i.e. $Prism_3$), Petersen graph, $3$-cube, $2$-truncated tetrahedron and a $10$-vertex graph (with $4$ vertices of degree five and all others of degree three), respectively.

\subsection{The cone $HCUT_6^3$}
\noindent There are two orbits of $4$-partitions of $V_6$ giving us the total of $65$ extreme rays.
There are $16$ orbits of facets; the representatives of corresponding orbits 
$F_i$, $1 \le i \le 16$ and their representation matrices are given in Table \ref{tab:tabl3P6}, where $\overline{ij}$ means the complement of a $2$-subset of $V_6$.

\begin{table}
\begin{center}
\tiny
\begin{tabular}{|c|ccccccccccccccc|c|c|c|}
\hline
&&&&&&&&&&&&&&&&&&\\[-1pt]
&$\overline{56}$&$\overline{46}$&$\overline{45}$&$\overline{36}$
&$\overline{35}$&$\overline{34}$&$\overline{26}$&$\overline{25}$
&$\overline{24}$&$\overline{23}$&$\overline{16}$&$\overline{15}$
&$\overline{14}$&$\overline{13}$&$\overline{12}$&Adj.&Size&Inc.\\
\hline
E1&0&0&0&0&0&0&0&0&1&1&0&0&1&1&0&58&45&993\\
E2&0&0&0&0&0&0&0&0&0&1&0&0&0&1&1&55&20&1113\\
\hline
F1&0&0&0&0&0&0&0&0&0&0&0&0&0&0&1&1526&15&49\\
F2&-1&0&1&0&1&0&0&1&0&0&0&1&0&0&0&703&30&41\\
F3&-1&0&1&0&1&0&1&0&1&1&2&1&0&0&-1&100&180&23\\
F4&-1&0&1&0&1&0&1&0&1&1&1&2&-1&1&0&37&360&19\\
F5&-1&-1&2&0&1&1&2&1&1&0&2&1&1&2&-2&31&360&18\\
F6&-1&-1&2&1&0&2&1&2&0&2&2&1&1&-1&-1&30&180&18\\
F7&-1&-1&2&1&0&2&2&1&1&-1&2&3&-1&1&0&23&360&16\\
F8&-1&-2&3&1&2&1&2&1&2&-1&2&1&2&3&-2&23&360&17\\
F9&-1&-2&3&2&1&2&2&1&2&-2&3&4&-1&1&1&23&180&15\\
F10&-1&-1&2&-1&2&2&2&1&1&1&2&1&1&1&-2&22&60&18\\
F11&-1&-1&2&0&1&1&1&2&2&-1&2&1&1&2&-1&18&360&16\\
F12&-1&-1&2&1&1&1&1&1&1&-1&1&1&1&-1&2&18&90&16\\
F13&-1&-1&2&1&0&2&1&2&0&0&2&1&1&-1&1&14&720&14\\
F14&-1&0&1&0&1&2&2&1&0&2&2&1&2&0&-2&14&360&14\\
F15&-1&-1&2&1&0&2&2&1&1&1&2&1&1&1&-2&14&360&14\\
F16&-1&0&1&0&1&0&1&0&1&1&1&0&1&1&0&14&90&14\\
\hline
\end{tabular}
 \begin{tabular}{|c|p{8pt}p{8pt}p{8pt}p{8pt}p{8pt}p{8pt}p{8pt}p{8pt}p{8pt}p{8pt}p{8pt}p{8pt}p{8pt}p{8pt}p{8pt}p{8pt}|c|c|} 
\hline
&F1&F2&F3&F4&F5&F6&F7&F8&F9&F10&F11&F12&F13&F14&F15&F16&Adj.&Size\\
\hline
F1& 14& 28& 144& 144& 192& 96& 120& 168& 36& 20& 120& 36& 240& 72& 72& 24&
1526&15\\
F2& 14& 25& 72& 96& 60& 24& 72& 48& 36& 10& 48& 12& 72& 48& 48& 18&703&30\\
F3& 12& 12& 14& 10& 10& 6& 10& 0& 0& 2& 2& 0& 12& 4& 4& 2&100&180\\
F4& 6& 8& 5& 2& 2& 2& 2& 2& 2& 0& 2& 0& 0& 2& 2& 0&37&360\\
F5& 8& 5& 5& 2& 3& 0& 0& 2& 0& 1& 1& 0& 0& 2& 2& 0&31&360\\
F6& 8& 4& 6& 4& 0& 2& 4& 2& 0& 0& 0& 0& 0& 0& 0& 0&30&180\\
F7& 5& 6& 5& 2& 0& 2& 0& 0& 1& 0& 0& 0& 2& 0& 0& 0&23&360\\
F8& 7& 4& 0& 2& 2& 1& 0& 0& 2& 0& 2& 1& 2& 0& 0& 0&23&360\\
F9& 3& 6& 0& 4& 0& 0& 2& 4& 0& 0& 0& 0& 0& 2& 2& 0&23&180\\
F10& 5& 5& 6& 0& 6& 0& 0& 0& 0& 0& 0& 0& 0& 0& 0& 0&22&60\\
F11& 5& 4& 1& 2& 1& 0& 0& 2& 0& 0& 0& 1& 2& 0& 0& 0&18&360\\
F12& 6& 4& 0& 0& 0& 0& 0& 4& 0& 0& 4& 0& 0& 0& 0& 0&18&90\\
F13& 5& 3& 3& 0& 0& 0& 1& 1& 0& 0& 1& 0& 0& 0& 0& 0&14&720\\
F14& 3& 4& 2& 2& 2& 0& 0& 0& 1& 0& 0& 0& 0& 0& 0& 0&14&360\\
F15& 3& 4& 2& 2& 2& 0& 0& 0& 1& 0& 0& 0& 0& 0& 0& 0&14&360\\
F16& 4& 6& 4& 0& 0& 0& 0& 0& 0& 0& 0& 0& 0& 0& 0& 0&14&90\\
\hline
Size&15&30&180&360&360&180&360&360&180&60&360&90&720&360&360&90&&4065\\
\hline
\end{tabular}
\caption{Representatives of orbits of extreme rays and facets of $HCUT^3_6$, followed by the representation matrix of the ridge graph of $HCUT_6^3$}
\label{tab:tabl3P6}
\end{center}
\end{table}

The orbit $F_1$ of the non-negativity facets form a dominating
clique in the ridge graph of the cone $HCUT^3_6$.
Also, the ridge graph of the cone $HCUT^3_6$, restricted on the orbits
$F_1$ and $F_2$, coincides with the ridge graph of the cone $HMET^3_6$.
The complement of the {\em local graph} (i.e. of the graph induced by all neighbors of an representative of the orbit) for the orbits $F_{16}, F_{15}, F_{14}, F_{13}, F_{12}$
of small adjacency have, respectively $8, 8, 9, 7, 12$ vertices. It
is $4K_2$ for $F_{16}$ and some connected planar graphs for $F_{15}, F_{14}, F_{13}$;
for $F_{12}$ it is $K_8-C_8$ on the set $V_8$ plus four 
pendent edges, which are incident with vertices $1, 3, 5, 7$, respectively.

\section{The case $n=m+3$}\label{case-m-m+3}
The dimension of the super-metric cone $SMET^{m,s}_{n=m+3}$ is ${n\choose{n-2}}=\frac{n(n-1)}{2}$, i.e. the same as of $MET_n$.

This correspondence allow us to replace the graph $G_v$ by a simpler graph: any facet or extreme ray of the cone $SMET^{m, s}_{n}$ is given by a vector, say, $v$, indexed by $(m+1)$-subsets of $V_{n}$. It can be seen also as a function on $2$-subsets of $V_{n}$, which are complements of $(m+1)$-subsets of $V_{n=m+3}$. So, to every $0,1$ valued extreme ray of $SMET^{m, s}_{n}$, one can associate a set of pairs $(ij)$ and this set of pairs is edge-set of a graph $H_v$, such that the graph $G_v$ is the line graph of $H_v$. If some vertices are isolated, then we remove them.

For example, if $v$ is a cut $\delta_{\{1, 2, 3\}, \{4\}}$, then $v$ has the support $\{\{1,4\}, \{2,4\}, \{3,4\}\}$, 
i.e. the complements are $\{\{2,3\}, \{1,3\}, \{1,2\}\}$; so, $H_v$ is the complete graph on $\{1, 2, 3\}$.

If $v=\alpha(\{1, 2, 3\}, \{4\}, \{5\})$ (an extreme ray of $HMET^2_5$), then its support is $\{\{1,4,5\}, \{2,4,5\}, \{3,4,5\}\}$, i.e. the complements are $\{\{2,3\}, \{1,3\}, \{1,2\}\}$; so, $H_v$ is the complete graph on vertices $\{1, 2, 3\}$ again. In fact, $\alpha(\{1, 2, 3\}, \{4\}, \{5\})$ is zero-extension of $\delta_{\{1, 2, 3\}, \{4\}}$.

For any $0,1$ valued vector $d$, the graph $H_{d^{ze}}$ of the zero-extension of $d$, is equal to $H_{d}$.

In these terms, we have
\begin{enumerate}
\item All extreme rays of $SMET^{1, 1}_4=MET_4$ have $H_v=K_3$ or $K_{2,2}(=C_4)$.
\item All extreme rays of $HMET^{2}_5$ have $H_v=K_3(=C_3)$, $K_{2,2}(=C_4)$; $\overline{C_5}(=C_5)$.
\item All $0,1$ extreme rays of $SMET^{2, 2}_5$ have $H_v=K_4$ or $K_{2,2,1}$.
\item All $0,1$ extreme rays of $HMET^{3}_6$ have $H_v=K_3(=C_3)$, $K_{2,2}(=C_4)$; $\overline{C_5}(=C_5)$; $C_6$.
\item All $0,1$ extreme rays of $SMET^{3,2}_6$ have $H_v=K_4$, $K_{2,2,1}$; $\overline{C_6}$(=$Prism_3$), $\overline{C_1+C_5}$, $\overline{C_3+C_3}$.
\item All $0,1$ extreme rays of $SMET^{3,3}_6$ have $H_v=K_5$ or $K_{2,2,2}$.
\item All $0,1$ extreme rays of $HMET^{4}_7$ have $H_v=K_3(=C_3)$, $K_{2,2}(=C_4)$; $C_5$, $C_6$, $C_7$.
\item All $0,1$ extreme rays of $SMET^{4,2}_7$ have $H_v=K_4$, $K_{2,2,1}$; $\overline{C_6}$, $\overline{C_1+ C_5}$, $\overline{C_3+C_3}$; $\nabla C_6$ or the graphs
\begin{center}
\epsfxsize=90mm
\epsfbox{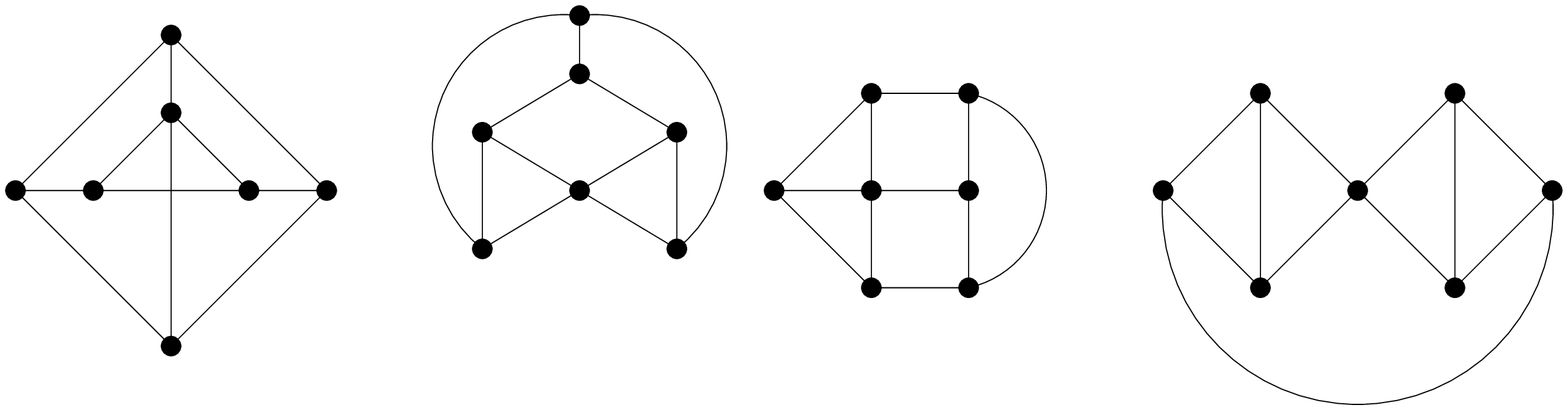}
\end{center}
\item All $0,1$ extreme rays of $SMET^{4,3}_7$ have $H_v=K_5$, $K_{2,2,2}$; $\overline{C_7}$, $\overline{C_1+C_6}$, $\overline{C_2+C_5}$, $\overline{C_3+C_4}$, $\overline{C_1+C_3+C_3}$.
\item All $0,1$ extreme rays of $SMET^{4, 4}_7$ have $H_v=K_6$ or $K_{2,2,2,1}$.
\item All $0,1$ extreme rays of $HMET^{5}_8$ have $H_v=K_3(=C_3)$, $K_{2,2}(=C_4)$; $C_5$, $C_6$, $C_7$, $C_8$.
\item Some $0,1$ extreme rays of $SMET^{5,2}_8$ have $H_v=K_4$, $K_{2,2,1}$; $\overline{C_6}$, $\overline{C_1+C_5}$, $\overline{C_3+C_3}$, $\nabla C_7$; four graphs depicted above for $SMET^{4,2}_7$, $\nabla C_6$, $3$-cube, or the graphs
\begin{center}
\epsfxsize=120mm
\epsfbox{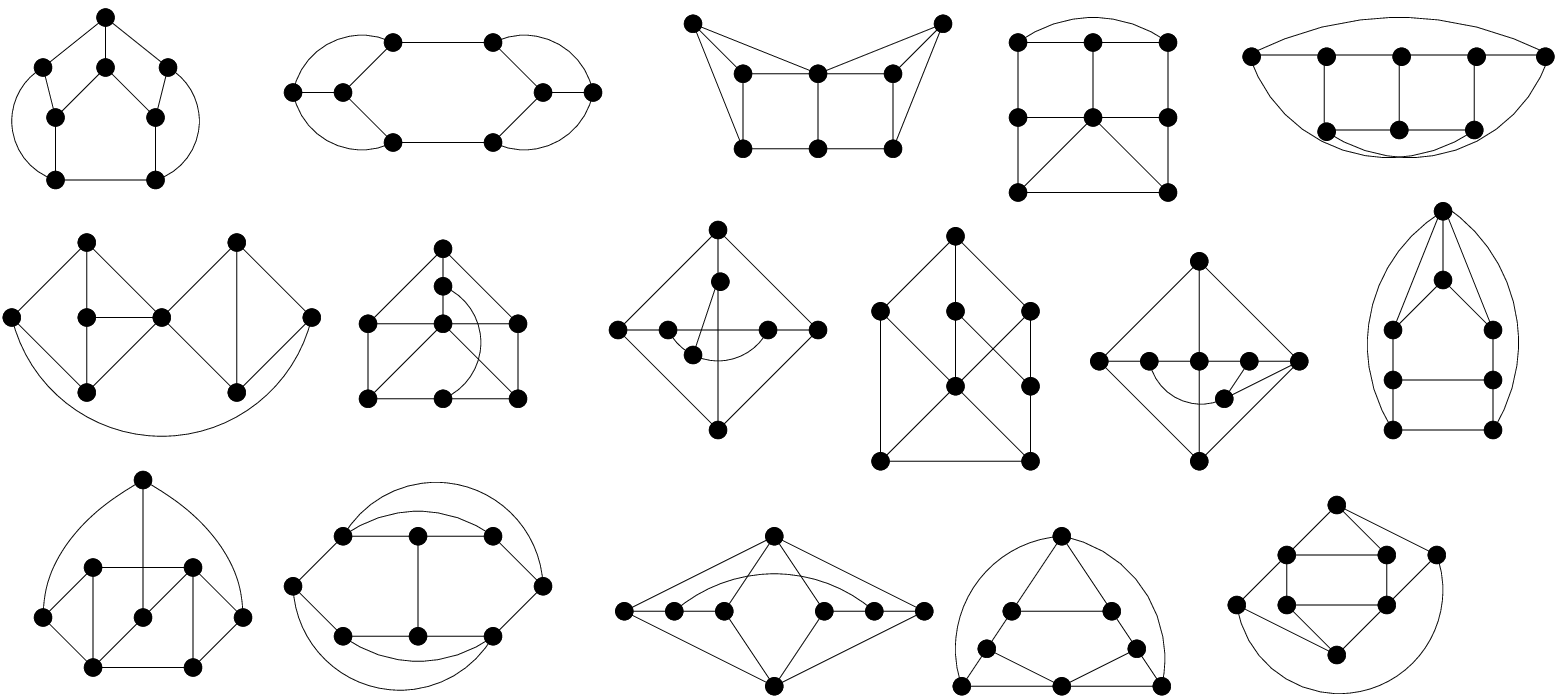}
\end{center}
\item Some $0,1$ extreme rays of $SMET^{5,3}_8$ have $H_v=K_5$, $K_{2,2,2}$; $K_{4,4}$, $\overline{C_7}$, $\overline{C_1+C_6}$, $\overline{C_2+C_5}$, $\overline{C_3+C_4}$, $\overline{C_1+C_3+C_3}$; complement of $3$-cube or the complement of the graphs
\begin{center}
\epsfxsize=120mm
\epsfbox{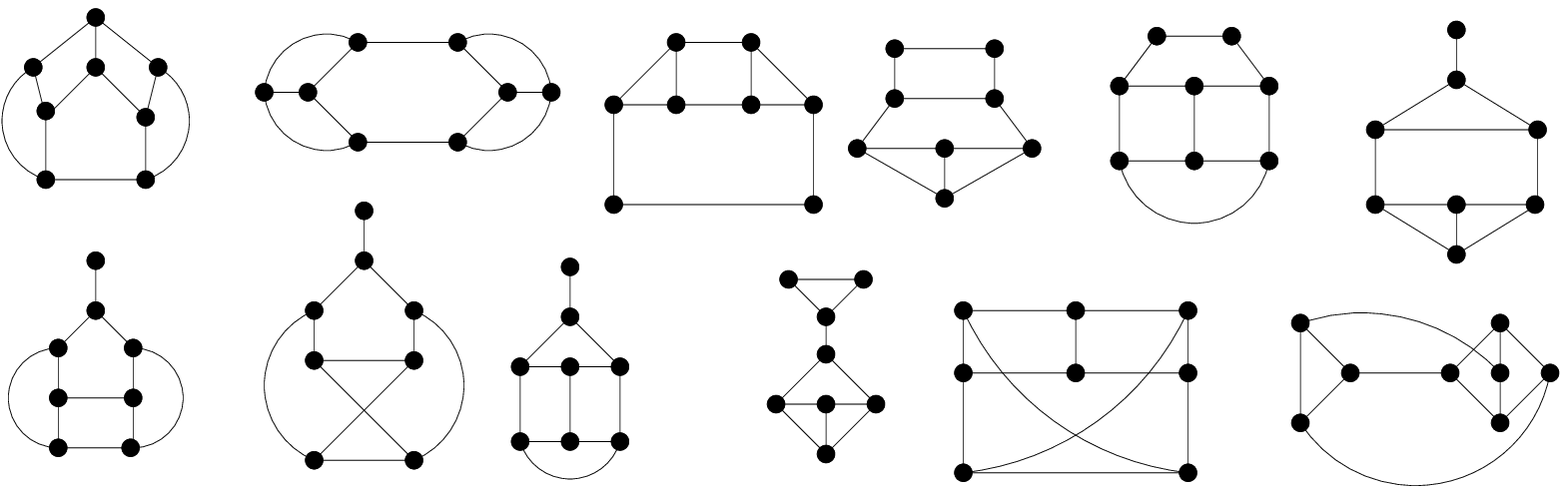}
\end{center}
\item Some $0,1$ extreme rays of $SMET^{5, 4}_8$ have $H_v=K_6$, $K_{2,2,2,1}$; $\overline{C_8}$, $\overline{C_1+C_7}$, $\overline{C_2+C_6}$, $\overline{C_3+C_5}$, $\overline{C_4+C_4}$, $\overline{C_1+C_3+C_4}$, $\overline{C_2+C_3+C_3}$.
\item All $0,1$ extreme rays of $SMET^{5, 5}_8$ have $H_v=K_7$ or $K_{2,2,2,2}$.
\end{enumerate}
%Note that for cones $SMET^{5,2}_8$, $SMET^{5,3}_8$, $SMET^{5,4}_8$ we may not have the complete list of extreme rays.

\begin{conj}
In terms of graph $H_v$, associated to $0,1$ vectors, we have for the cone $SMET^{m, s}_{n=m+3}$:

(i) For $s=m$, there are two orbits of $0,1$ extreme rays: $K_{n-1}$ and $K_n-\lfloor\frac{n}{2}\rfloor K_2$.

(ii) for $s=m-1$, besides zero-extension from $SMET^{m-1, m-1}_{n-1}$, all 
$0,1$ valued extreme rays for $m=2, 3, 4$ and some for $m=5$ have $H_v$ being the complement of the union of disjoint circuits with lengths partitioning $n=m+3$, but lengths-vectors (the lengths are not increasing) $(\dots,1,1)$, $(\dots,2,1)$, $(\dots,2,2)$ and
$(n\leq 4)$, $(n-1\leq 4,1)$, $(n-2\leq 4,2)$ are excluded.

%(iii) For $s=1$ and $3\leq i\leq m+3$, there is an orbit of extreme ray with $H_v=C_i$.

%(iv) For $s=m-2$ or $2$: besides zero-extension from (ii) for $m-1$, there is an orbit with $H_v=K_n-\lfloor n/2\rfloor K_2-2K_{\lfloor n/2\rfloor}$ (i.e. $C_6$, $\nabla C_6$, $3$-cube, \dots for $n=6,7,8$) starting from $n=6$ or $7$ for $s=m-2$ or $2$, respectively.
\end{conj}
We checked this conjecture up to $n=12$. 
%Let us call zero-extensions of the rays in (i) above simplex and hyperoctahedral ray
%A zero-extension of the extreme ray $v$ with $H_v=K_1+K_{n-1}$ is an extreme ray of $SMET^{m,m+1}_{}$
%call it simplex 

\begin{conj}
The complete list of extreme rays of $HMET^{m}_{n=m+3}$ consists of

(i) $0,1$ valued extreme rays with $H_v=C_i$ for $3\leq i\leq m+3$ and 

(ii) $0,1,2$ valued extreme rays $v$ represented by
$H_v=C_{1,2,\dots,i} + P_{i,i+1,\dots,k} + C_{k,k+1,\dots,j}$
(with $3\leq i<k$, $j\geq k+2$, and $i+k-1\leq j\leq n$), having value $2$ on 
edges of the path (so, besides zero-extensions, they are those with $j=n$, 
i.e. with all $3\leq i<k\leq n-i+1$).
\end{conj}

Another interesting fact is that the found minimal incidence number for extreme rays of $SMET^{m,2}_{m+3}$ is ${{m+3}\choose 2}$ and not (apriori possible) minimum ${{m+3}\choose 2}-1$, which occur for other super-metric cones.

%\begin{conj}
%The $m$-vectors of $m+3$ sets have $H_v=\overline{3C_1}, \overline{2C_2}$. All other extreme rays of $HMET^{m}_{m+3}$ have $H_v=C_i$ for all $3\leq i\leq m+3$.
%\end{conj}
%This conjecture was checked for $MET_4$, $HMET^{2}_5$, $HMET^{3}_6$, $HMET^{4}_7$, $HMET^{5}_8$

\begin{conj}
The number of extreme rays, the number of orbits, and the diameter of $SMET^{2, s}_5$ are:

{\scriptsize
\begin{center}
\begin{tabular}{|l|l|l|l|}
\hline
$(1170, 16, 5)$ if $s\in ]0,1[$&$(37, 3, 2)$ ~~~~if $s=1$&$(2462, 35, 5)$ if $s\in ]1, \frac{3}{2}[$&$(1442, 25, 4)$ if $s=\frac{3}{2}$\\
\hline
$(2102, 31, 5)$ if $s\in ]\frac{3}{2}, \frac{5}{3}[$&$(1742, 28, 5)$ if $s=\frac{5}{3}$&$(1862, 30, 5)$ if $s\in ]\frac{5}{3}, 2[$&$(132, 6, 2)$ ~~~if $s=2$\\
\hline
\end{tabular}
\end{center}
}

\end{conj}
This conjecture was verified for about one hundred values of $s$.

The representatives of orbits of facets and extreme rays of the cones $SMET^{2,2}_5$, $SCUT^{2,2}_5$, $SMET^{3,3}_6$, $SCUT^{3,3}_6$, $SMET^{4,4}_7$, $SMET^{3,2}_6$ are presented in Tables \ref{table22smet5}, \ref{tab:22SCUT5}, \ref{table33smet6}, \ref{tab:33SCUT6}, \ref{tabl44smet7}, \ref{tabl32smet6}, respectively.

\begin{table}
\begin{center}
\tiny
\begin{tabular}{|c|cccccccccc|c|c|c|}
\hline
&&&&&&&&&&&&&\\[-1pt]
&$\overline{45}$&$\overline{35}$&$\overline{34}$&$\overline{25}$
&$\overline{24}$&$\overline{23}$&$\overline{15}$&$\overline{14}$
&$\overline{13}$&$\overline{12}$&Adj.&Size&Inc.\\
\hline
F1&-2&0&1&0&1&0&0&1&0&0&19&20&70\\
\hline
E1&0&0&1&0&1&1&0&1&1&1&92&5&16\\
E2&0&2&2&2&2&1&2&2&1&1&32&10&12\\
E3&0&1&1&1&1&0&1&1&1&1&28&15&12\\
E4&1&1&2&2&1&2&2&2&1&1&25&12&10\\
E5&1&2&4&3&3&2&3&4&4&1&13&60&10\\
E6&0&4&4&4&4&2&4&4&2&5&13&30&10\\
\hline
\end{tabular}
\begin{tabular}{|c|cccccc|c|c|}
\hline
&E1&E2&E3&E4&E5&E6&Adj.&Size\\
\hline
E1& 4& 10& 12& 12& 36& 18&92&5\\
E2& 5& 0& 3& 6& 12& 6&32&10\\
E3& 4& 2& 2& 4& 12& 4&28&15\\
E4& 5& 5& 5& 0& 5& 5&25&12\\
E5& 3& 2& 3& 1& 2& 2&13&60\\
E6& 3& 2& 2& 2& 4& 0&13&30\\
\hline
Size&5&10&15&12&60&30&&132\\
\hline
\end{tabular}
\caption{Representatives of orbits of facets and extreme rays of $SMET^{2,2}_5$, followed by the representation matrix of skeleton graph of $SMET^{2,2}_5$}\label{table22smet5}
\end{center}
\end{table}

\begin{table}
\begin{center}
\tiny
\begin{tabular}{|c|cccccccccc|c|c|c|}
\hline
&&&&&&&&&&&&&\\[-1pt]
&$\overline{45}$&$\overline{35}$&$\overline{34}$&$\overline{25}$
&$\overline{24}$&$\overline{23}$&$\overline{15}$&$\overline{14}$
&$\overline{13}$&$\overline{12}$&Adj.&Size&Inc.\\
\hline
$E_1$&0&1&1&1&1&0&1&1&1&1&19&15&104\\
$E_2$&0&0&1&0&1&1&0&1&1&1&19&5&134\\
\hline
$F_1$&-2&0&1&0&1&0&0&1&0&0&61&20&13\\
$F_2$&-1&-3&1&1&1&1&1&1&1&-1&22&60&11\\
$F_3$&-3&-3&3&-3&3&3&1&1&1&1&15&20&10\\
$F_4$&-1&-1&1&1&1&1&1&1&1&-3&13&30&10\\
$F_5$&-3&-3&3&1&1&1&1&3&3&-3&9&60&9\\
$F_6$&-1&-1&3&5&3&3&5&3&3&-15&9&30&9\\
\hline
\end{tabular}
\begin{tabular}{|c|cccccc|c|c|}
\hline
&$F_1$&$F_2$&$F_3$&$F_4$&$F_5$&$F_6$&Adj.&Size\\
\hline
$F_1$& 13& 21& 6& 6& 9& 6&61&20\\
$F_2$& 7& 6& 2& 3& 2& 2&22&60\\
$F_3$& 6& 6& 0& 0& 3& 0&15&20\\
$F_4$& 4& 6& 0& 0& 2& 1&13&30\\
$F_5$& 3& 2& 1& 1& 2& 0&9&60\\
$F_6$& 4& 4& 0& 1& 0& 0&9&30\\
\hline
Size&20&60&20&30&60&30&&220\\
\hline
\end{tabular}
\caption{Representatives of orbits of extreme rays and facets of $SCUT^{2,2}_5$, followed by the representation matrix of the ridge graph of $SCUT^{2,2}_5$}\label{tab:22SCUT5}
\end{center}
\end{table}

\begin{table}
\begin{center}
\tiny
\begin{tabular}{|c|ccccccccccccccc|c|c|c|}
\hline
&&&&&&&&&&&&&&&&&&\\[-1pt]
&$\overline{56}$&$\overline{46}$&$\overline{45}$&$\overline{36}$
&$\overline{35}$&$\overline{34}$&$\overline{26}$&$\overline{25}$
&$\overline{24}$&$\overline{23}$&$\overline{16}$&$\overline{15}$
&$\overline{14}$&$\overline{13}$&$\overline{12}$&Adj.&Size&Inc.\\
\hline
F1&-3&0&1&0&1&0&0&1&0&0&0&1&0&0&0&29&30&594\\
\hline
E1&0&0&1&0&1&1&0&1&1&1&0&1&1&1&1&650&6&25\\
E2&0&1&1&1&1&0&1&1&1&1&1&1&1&1&0&449&15&24\\
E3&0&3&3&3&3&0&3&3&3&3&3&3&3&3&4&93&45&18\\
E4&1&1&1&2&2&2&2&2&2&1&2&2&2&1&1&57&10&18\\
E5&0&2&2&2&2&1&2&2&1&1&2&2&2&2&2&56&60&17\\
E6&1&1&3&2&2&2&2&3&3&1&2&3&3&1&3&51&90&18\\
E7&1&1&2&2&1&2&2&2&1&1&2&2&2&2&2&30&72&15\\
E8&0&3&3&3&3&2&3&3&2&4&3&3&2&4&4&27&60&16\\
E9&0&4&4&4&4&2&4&4&2&5&4&4&4&5&5&23&180&15\\
E10&1&2&4&3&3&2&3&4&4&1&3&4&4&3&4&18&360&15\\
E11&1&3&3&4&4&4&4&4&5&2&4&4&5&2&5&18&180&15\\
E12&3&3&3&6&6&6&6&6&6&3&6&6&6&3&7&14&60&14\\
\hline
\end{tabular}
\begin{tabular}{|c|cccccccccccc|c|c|}
\hline
&E1&E2&E3&E4&E5&E6&E7&E8&E9&E10&E11&E12&Adj.&Size\\
\hline
E1& 5& 15& 30& 10& 50& 90& 60& 30& 90& 180& 60& 30&650&6\\
E2& 6& 14& 39& 6& 36& 48& 24& 36& 72& 72& 72& 24&449&15\\
E3& 4& 13& 2& 2& 4& 16& 8& 4& 8& 24& 8& 0&93&45\\
E4& 6& 9& 9& 0& 6& 9& 0& 0& 0& 0& 18& 0&57&10\\
E5& 5& 9& 3& 1& 0& 6& 6& 0& 6& 12& 6& 2&56&60\\
E6& 6& 8& 8& 1& 4& 4& 0& 4& 4& 8& 4& 0&51&90\\
E7& 5& 5& 5& 0& 5& 0& 0& 0& 5& 5& 0& 0&30&72\\
E8& 3& 9& 3& 0& 0& 6& 0& 0& 3& 0& 3& 0&27&60\\
E9& 3& 6& 2& 0& 2& 2& 2& 1& 0& 4& 1& 0&23&180\\
E10& 3& 3& 3& 0& 2& 2& 1& 0& 2& 2& 0& 0&18&360\\
E11& 2& 6& 2& 1& 2& 2& 0& 1& 1& 0& 0& 1&18&180\\
E12& 3& 6& 0& 0& 2& 0& 0& 0& 0& 0& 3& 0&14&60\\
\hline
Size&6&15&45&10&60&90&72&60&180&360&180&60&&1138\\
\hline
\end{tabular}
\caption{Representatives of orbits of facets and extreme rays of $SMET^{3,3}_6$, followed by the representation matrix of the skeleton graph of $SMET^{3,3}_6$}\label{table33smet6}
\end{center}
\end{table}

\begin{table}
\begin{center}
\tiny
\begin{tabular}{|c|ccccccccccccccc|c|c|c|}
\hline
&&&&&&&&&&&&&&&&&&\\[-1pt]
&$\overline{56}$&$\overline{46}$&$\overline{45}$&$\overline{36}$
&$\overline{35}$&$\overline{34}$&$\overline{26}$&$\overline{25}$
&$\overline{24}$&$\overline{23}$&$\overline{16}$&$\overline{15}$
&$\overline{14}$&$\overline{13}$&$\overline{12}$&Adj.&Size&Inc.\\
\hline
$E_1$&0&0&1&0&1&1&0&1&1&1&0&1&1&1&1&20&6&125\\
$E_2$&0&1&1&1&1&0&1&1&1&1&1&1&1&1&0&20&15&108\\
\hline
$F_1$&-1&-1&0&-1&0&0&0&1&1&1&0&1&1&1&-2&25&60&17\\
$F_2$&-3&0&1&0&1&0&0&1&0&0&0&1&0&0&0&25&30&17\\
$F_3$&-1&-1&1&-1&1&1&0&0&0&0&0&0&0&0&1&14&60&14\\
\hline
\end{tabular}
\begin{tabular}{|c|ccc|c|c|}
\hline
&$F_1$&$F_2$&$F_3$&Adj.&Size\\
\hline
$F_1$& 12& 7& 6&25&60\\
$F_2$& 14& 5& 6&25&30\\
$F_3$& 6& 3& 5&14&60\\
\hline
Size&60&30&60&&150\\
\hline
\end{tabular}
\caption{Representatives of orbits of extreme rays and facets of $SCUT^{3,3}_6$, followed by the representation matrix of the ridge graph of $SCUT^{3,3}_6$}\label{tab:33SCUT6}
\end{center}
\end{table}

\begin{table}
\begin{center}
\tiny
\begin{tabular}{|c|ccccccccccccccc|c|c|c|}
\hline
&&&&&&&&&&&&&&&&&&\\[-1pt]
&$\overline{56}$&$\overline{46}$&$\overline{45}$&$\overline{36}$
&$\overline{35}$&$\overline{34}$&$\overline{26}$&$\overline{25}$
&$\overline{24}$&$\overline{23}$&$\overline{16}$&$\overline{15}$
&$\overline{14}$&$\overline{13}$&$\overline{12}$&Adj.&Size&Inc.\\
\hline
F1&-2&0&1&0&1&0&0&1&0&0&0&1&0&0&0&43&30&5404\\
F2&0&0&0&0&0&0&0&0&0&0&0&0&0&0&1&34&15&3195\\
\hline
E1&0&0&0&0&0&1&0&0&1&1&0&0&1&1&1&1642&15&31\\
E2&0&0&1&1&0&1&1&1&0&1&1&1&1&0&0&953&60&24\\
E3&0&0&0&1&1&1&1&1&1&0&1&1&1&0&0&696&10&24\\
E4&0&0&0&0&1&1&0&1&1&0&0&1&1&1&1&274&90&24\\
E5&0&0&1&2&0&0&2&1&1&2&2&1&1&2&0&248&90&21\\
E6&0&0&0&0&2&2&0&2&2&1&0&2&2&1&1&183&60&23\\
E7&0&0&1&1&0&1&1&1&0&0&1&1&1&1&1&125&72&20\\
E8&0&0&0&2&2&2&2&2&2&0&2&2&2&0&3&103&60&19\\
E9&0&0&2&2&0&2&2&2&0&2&2&2&2&0&3&92&360&19\\
E10&0&0&1&1&0&1&1&1&2&2&1&1&2&2&0&84&180&18\\
E11&0&0&1&2&0&3&2&1&2&1&2&1&3&3&0&73&720&18\\
E12&0&0&0&0&4&4&0&4&4&2&0&4&4&2&5&59&180&21\\
E13&0&0&1&3&0&2&3&1&3&2&3&1&3&2&0&46&360&17\\
E14&0&1&2&1&2&1&2&1&0&2&2&1&2&0&1&39&60&15\\
E15&0&0&1&0&1&2&0&2&1&2&0&2&2&1&1&35&72&20\\
E16&0&0&3&1&1&1&1&2&3&3&1&3&2&3&0&32&360&16\\
E17&0&1&2&3&4&1&4&2&0&4&4&4&2&0&2&30&360&15\\
E18&0&0&1&0&2&4&0&3&3&2&0&3&4&4&1&23&360&20\\
E19&0&1&1&4&4&2&5&5&0&2&5&5&2&0&3&23&360&15\\
E20&0&1&2&1&2&3&2&1&0&3&2&1&3&0&3&23&180&15\\
E21&0&1&1&2&2&0&3&3&2&1&3&3&2&1&0&23&90&15\\
E22&0&0&1&2&0&1&2&1&2&1&2&1&2&2&0&22&720&16\\
E23&0&0&2&6&0&4&6&2&6&4&6&2&6&7&0&22&720&16\\
E24&0&0&1&2&0&2&2&1&2&1&2&1&2&1&0&22&360&16\\
E25&0&0&2&2&0&3&2&2&1&1&2&2&3&3&0&22&360&16\\
E26&0&0&2&3&0&2&3&2&1&3&3&2&1&3&0&22&360&16\\
E27&0&0&2&3&2&5&3&2&5&5&3&3&3&0&0&22&360&17\\
E28&0&0&2&2&0&1&2&2&1&1&2&2&2&2&0&22&360&16\\
E29&0&0&2&6&0&7&6&2&6&4&6&2&6&4&0&22&360&16\\
E30&0&0&2&2&0&5&2&2&4&4&2&2&4&4&0&22&180&16\\
E31&0&0&5&1&1&2&1&4&3&4&1&5&5&4&0&20&720&16\\
E32&0&0&1&1&1&3&1&2&2&1&1&2&3&3&0&20&720&16\\
E33&0&0&3&3&1&5&3&2&5&5&3&3&2&1&0&20&720&16\\
E34&0&0&5&3&2&3&3&3&2&0&3&5&5&1&1&20&360&16\\
E35&0&0&1&2&1&3&2&2&2&0&2&2&3&3&0&18&720&16\\
E36&0&0&2&2&1&3&2&1&3&3&2&2&2&0&0&18&360&16\\
E37&0&0&1&2&1&1&2&2&2&0&2&2&2&2&0&18&360&16\\
E38&0&1&2&5&5&2&6&3&0&3&6&5&1&0&6&18&360&15\\
E39&0&1&4&3&6&3&4&2&4&0&4&6&0&6&2&18&360&15\\
E40&0&1&2&1&2&1&2&4&0&2&2&4&2&0&4&18&120&15\\
\hline
\end{tabular}
\caption{Representatives of orbits of facets and extreme rays of $SMET^{3,2}_6$}\label{tabl32smet6}
\end{center}
\end{table}

\begin{table}
\begin{center}
\tiny
\begin{tabular}{|c|p{7pt}p{3pt}p{3pt}p{3pt}p{3pt}p{3pt}p{3pt}p{3pt}p{3pt}p{3pt}p{3pt}p{3pt}p{3pt}p{3pt}p{3pt}p{3pt}p{3pt}p{3pt}p{3pt}p{3pt}p{3pt}|c|c|c|}
\hline
&&&&&&&&&&&&&&&&&&&&&&&&\\[-1mm]
&$\overline{67}$&$\overline{57}$&$\overline{56}$&$\overline{47}$&$\overline{46}$&$\overline{45}$&$\overline{37}$&$\overline{36}$&$\overline{35}$&$\overline{34}$&$\overline{27}$&$\overline{26}$&$\overline{25}$&$\overline{24}$&$\overline{23}$&$\overline{17}$&$\overline{16}$&$\overline{15}$&$\overline{14}$&$\overline{13}$&$\overline{12}$&Adj.&Size&Inc.\\
\hline
F1&-4&0&1&0&1&0&0&1&0&0&0&1&0&0&0&0&1&0&0&0&0&41&42&21363\\
\hline
E1&0&0&1&0&1&1&0&1&1&1&0&1&1&1&1&0&1&1&1&1&1&17163&7&36\\
E2&0&1&1&1&1&0&1&1&1&1&1&1&1&1&0&1&1&1&1&1&1&1486&105&30\\
E3&0&3&3&3&3&0&3&3&3&3&3&3&3&3&4&3&3&3&3&4&4&1314&105&26\\
E4&0&2&2&2&2&0&2&2&2&2&2&2&2&2&1&2&2&2&2&1&1&1228&105&32\\
E5&0&2&2&2&2&1&2&2&1&2&2&2&2&1&2&2&2&2&2&1&1&343&252&30\\
E6&0&8&8&8&8&0&8&8&8&8&8&8&8&8&4&8&8&8&8&4&9&294&315&26\\
E7&0&4&4&4&4&2&4&4&2&2&4&4&4&4&4&4&4&4&4&4&5&238&210&24\\
E8&0&4&4&4&4&2&4&4&2&5&4&4&4&5&5&4&4&4&5&5&3&153&630&24\\
E9&1&2&4&3&3&2&3&4&4&1&3&4&4&3&4&3&4&4&3&4&2&120&1260&26\\
E10&1&1&2&2&1&2&2&2&1&2&2&2&2&1&2&2&2&2&2&1&1&112&360&28\\
E11&2&2&4&4&2&5&4&4&4&4&4&4&5&5&2&4&4&5&5&2&5&111&1260&24\\
E12&1&1&3&2&2&2&2&3&3&1&2&3&3&1&3&2&3&3&2&3&3&108&630&24\\
E13&0&8&8&8&8&4&8&8&4&8&8&8&8&4&8&8&8&8&8&4&9&95&1260&24\\
E14&0&8&8&8&8&10&8&8&10&10&8&8&10&10&11&8&8&4&4&8&8&93&630&22\\
E15&0&9&9&9&9&10&9&9&10&10&9&9&3&6&9&9&9&9&6&3&9&70&1260&24\\
E16&1&1&1&2&2&2&2&2&2&1&2&2&2&1&1&2&2&2&2&2&2&63&70&24\\
E17&2&2&4&4&2&4&4&4&2&2&4&4&4&4&4&4&4&4&4&4&5&61&252&22\\
E18&2&4&8&6&6&4&6&8&8&2&6&8&8&6&8&6&8&8&6&8&9&49&1260&22\\
E19&0&8&8&8&8&11&8&8&11&11&8&8&11&11&11&8&8&6&6&6&6&47&105&24\\
E20&3&3&8&4&8&8&5&5&5&4&5&8&8&8&3&5&8&8&8&3&8&46&630&26\\
E21&2&2&6&4&4&4&4&6&6&2&4&6&6&3&6&4&6&6&3&6&5&44&1260&25\\
E22&1&2&4&3&3&3&3&4&3&3&3&4&4&1&4&3&4&4&2&3&4&40&2520&24\\
E23&0&3&3&3&3&2&3&3&2&4&3&3&2&4&4&3&3&3&4&4&4&40&420&22\\
E24&4&4&8&8&4&8&8&8&4&8&8&8&8&4&8&8&8&8&8&4&9&36&2520&22\\
E25&2&6&9&8&8&4&8&8&8&8&8&9&9&4&6&8&9&9&8&2&9&35&2520&24\\
E26&2&6&9&8&8&6&8&8&6&8&8&9&9&2&8&8&9&9&8&2&9&35&840&24\\
E27&2&4&8&6&6&6&6&8&9&4&6&8&9&4&9&6&8&9&4&9&9&31&840&23\\
E28&1&5&5&6&6&6&6&6&8&4&6&6&8&4&8&6&6&8&4&8&8&31&420&23\\
E29&10&10&10&3&6&9&9&6&9&9&9&9&3&9&9&9&9&9&9&3&6&28&5040&22\\
E30&12&12&12&12&12&8&12&16&16&4&4&16&16&12&16&8&17&16&12&16&16&24&5040&21\\
E31&12&12&12&12&15&3&12&9&9&15&6&16&12&15&15&6&16&12&15&15&16&24&2520&21\\
E32&1&3&3&4&4&4&4&4&5&2&4&4&5&2&5&4&4&5&4&5&5&24&1260&21\\
E33&12&12&12&12&12&3&12&6&9&9&6&14&12&12&12&6&14&12&12&12&14&24&1260&21\\
E34&12&12&25&16&16&16&16&24&24&8&16&24&24&8&24&8&24&24&16&24&24&24&1260&21\\
E35&12&12&28&12&28&28&18&18&18&18&18&27&27&27&9&18&27&27&27&9&27&24&420&21\\
E36&12&12&22&12&22&22&18&18&18&18&18&18&18&18&9&18&18&18&18&9&9&24&140&21\\
E37&3&3&3&6&6&6&6&6&6&3&6&6&6&3&7&6&6&6&6&7&7&20&420&20\\
\hline
\end{tabular}
\caption{Representatives of orbits of facets and extreme rays of $SMET^{4,4}_7$}\label{tabl44smet7}
\end{center}
\end{table}

\section{The case $m=s$ and $n\geq m+4$}
As we saw above, the case $n=m+2$ was completely solvable and the case $n=m+3$ was computationally reasonable. But already for $n=m+4$ the combinatorial explosion happens (see, for example, $HCUT^2_6$, $SMET^{2, 2}_6$, and $SCUT^{2,2}_6$).

A remarkable feature of the cone $SCUT^{1,1}_n=CUT_n$ is that its skeleton is a complete graph. We checked that this happens also for $SCUT^{m, s}_n$ with $(m,s,n)$=$(2, 2,5)$, $(3, 3, 6)$, $(4,4,7)$, $(2,2,6)$, $(2,2,7)$, but not for $(5,5,8)$, $(3,3,7)$ and $(3,2,6)$. Remark that $SCUT^{4,4}_7$ is the smallest example, when $SCUT^{m,s}_n$ has facets, {\it which are not adjacent} to any of $(m,s)$-simplex facets.

The cone $SMET^{m, m}_n$ shares with $SMET^{1,1}_n=MET_n$ the property that the non-negativity inequality is redundant.
We expect that it shares also the following conjectured property:
\begin{conj}
The cone $SMET^{m,m}_n$ has extreme rays without components zeros.
\end{conj}
This conjecture holds for the cones $MET_n=SMET^{1, 1}_n$, $SMET^{2, 2}_5$, $SMET^{3, 3}_6$, $SMET^{4, 4}_7$, $SMET^{2, 2}_6$, $SMET^{3, 3}_7$, $SMET^{2,2}_7$.

\begin{proposition}
The vector $v$, defined below, is an extreme ray of the cone $SMET^{m,m}_n$:
\begin{equation*}
v(S)=\left\lbrace\begin{array}{rcl}
1&\mbox{~if~}&1\in S,\\
0&\mbox{~if~}&1\notin S
\end{array}\right.
\end{equation*}
for any $(m+1)$-subset of $\{1, \dots, n\}$. Its graph $G_v$ is Johnson graph $J(n-1,m)$ and the orbit represented by $v$, has size $n$.
\end{proposition}
\proof If $T=\{x_1, \dots, x_{m+2}\}$ is an $(m+2)$-set, then $v$ is incident to the $(m,m)$-simplex facet $(s+1)d_{x_i}\leq \Sigma_T$ if and only if $1\notin T$ or $1\in T$ and $x_i\not= 1$. Let $d\in SMET^{m,m}_{n}$ be a ray incident to all those facets.

If $1\notin T$, then, by summing all $(m,m)$-simplex equalities with support in $T$, one obtains $(m+1)\Sigma_T=(m+2)\Sigma_T$, i.e. $\Sigma_T=0$. So, $d(S)=0$ if $1\notin S$.

If $1\in T$, we get $(m+1)d_{x}=\Sigma_T$ for all $x\in T-\{1\}$. Considering all $1\in T$, we obtain $d(S)=\lambda$ for all sets $S$ with $1\in S$. So, $d=\lambda v$, which proves the result.

\end{document}